\RequirePackage{fix-cm}
\documentclass[smallextended,envcountsect]{svjour3}       
\smartqed 

\usepackage{graphicx}
\usepackage{mathptmx} 

\usepackage{amsmath,amssymb}
\usepackage{array}
\usepackage{enumitem}
\usepackage{xcolor}
\usepackage{colortbl}
\usepackage{multirow}
\usepackage{bm}
\usepackage{subfigure}
\usepackage{algorithm}
\usepackage[noend]{algpseudocode}
\usepackage{algpseudocode}
\usepackage{setspace}
\usepackage{epstopdf}
\usepackage[normalem]{ulem}
\let\Algorithm\algorithm
\renewcommand\algorithm[1][]{\Algorithm[#1]\setstretch{1.0}}

\providecommand{\ie}		{\emph{i.e\@.}\xspace}
\providecommand{\eg}		{\emph{e.g\@.}\xspace}

\providecommand{\myurl}[1][]	{\texttt{web.eecs.umich.edu/$\sim$fessler#1}\xspace}

\providecommand{\onweb}[1]	{Available from \myurl.}

\long\def\comment#1{}

\providecommand{\bcent}		{\begin{center}}
\providecommand{\ecent}		{\end{center}}
\providecommand{\benum}		{\begin{enumerate}}
\providecommand{\eenum}		{\end{enumerate}}
\providecommand{\bitem}		{\begin{itemize}}
\providecommand{\eitem}		{\end{itemize}}
\providecommand{\bvers}		{\begin{verse}}
\providecommand{\evers}		{\end{verse}}
\providecommand{\btab}		{\begin{tabbing}}	
\providecommand{\etab}		{\end{tabbing}}

\newcounter{blist}
\providecommand{\blistmark}	{\makebox[0pt]{$\bullet$}}
\providecommand{\blistitemsep}	{0pt}
\providecommand{\blist}[1][]	{%
\begin{list}{\blistmark}{%
\usecounter{blist}%
\setlength{\itemsep}{\blistitemsep}%
\setlength{\parsep}{0pt}%
\setlength{\parskip}{0pt}%
\setlength{\partopsep}{0pt}%
\setlength{\topsep}{0pt}%
\setlength{\leftmargin}{1.2em}%
\setlength{\labelsep}{0.5\leftmargin}
\setlength{\labelwidth}{0em}%
#1}
}
\providecommand{\elist}		{\end{list}}

\providecommand{\blistitemsep}	{0pt}
\providecommand{\bjfenum}[1][]	{%
\begin{list}{\bcolor{\arabic{blist}.} }{%
\usecounter{blist}%
\setlength{\itemsep}{\blistitemsep}%
\setlength{\parsep}{0pt}%
\setlength{\parskip}{0pt}%
\setlength{\partopsep}{0pt}%
\setlength{\topsep}{0pt}%
\setlength{\leftmargin}{0.0em}%
\setlength{\labelsep}{1.0\leftmargin}
\setlength{\labelwidth}{0pt}%
#1}
}

\newcounter{blistAlph}
\providecommand{\blistAlph}[1][]
{\begin{list}{\makebox[0pt][l]{\Alph{blistAlph}.}}{%
\usecounter{blistAlph}%
\setlength{\itemsep}{0pt}\setlength{\parsep}{0pt}%
\setlength{\parskip}{0pt}\setlength{\partopsep}{0pt}%
\setlength{\topsep}{0pt}%
\setlength{\leftmargin}{1.2em}%
\setlength{\labelsep}{1.0\leftmargin}
\setlength{\labelwidth}{0.0\leftmargin}#1}%
}

\newcounter{blistRoman}
\providecommand{\blistRoman}[1][]
{\begin{list}{\Roman{blistRoman}.}{%
\usecounter{blistRoman}%
\setlength{\itemsep}{0.5em}\setlength{\parsep}{0pt}%
\setlength{\parskip}{0pt}\setlength{\partopsep}{0pt}%
\setlength{\topsep}{0pt}%
\setlength{\leftmargin}{4em}%
\setlength{\labelsep}{0.4\leftmargin}
\setlength{\labelwidth}{0.6\leftmargin}#1}%
}

\usepackage{bbm} 

\providecommand{\inprod}[2]	{\xmath{\mathop{\langle #1,\, #2 \rangle}\nolimits}}
\providecommand{\Inprod}[2]	{\xmath{\left\langle #1,\ #2 \right\rangle}}

\let\equivsave\equiv
\def\equiv{\xmath{\equivsave}}

\providecommand{\ba}[1]		{\left[ \begin{array}{#1}}
\providecommand{\ea}		{\end{array} \right]}
\providecommand{\be}		{\begin{equation}}
\providecommand{\ee}[1]		{\label{#1}\end{equation}}
\providecommand{\bea}		{\begin{eqnarray}}
\providecommand{\eea}[1]	{\label{#1}\end{eqnarray}}
\providecommand{\beas}		{\begin{eqnarray*}}
\providecommand{\eeas}		{\end{eqnarray*}}
\providecommand{\beals}[1][1]	{\begin{alignat*}{#1}}	
\providecommand{\eeals}		{\end{alignat*}}

\providecommand{\berr}[2]{
\bgroup
\renewcommand{\theequation}{#1}
\be
#2
\ee{e,#1}
\egroup
\ignorespaces
}

\providecommand{\bearr}[2]{
\bgroup
\renewcommand{\theequation}{#1}
\bea
#2
\eea{e,#1}
\egroup
\ignorespaces
}

\providecommand{\inmath}	{\ensuremath}
\providecommand{\xmath}[1]	{\inmath{#1}\xspace}
\providecommand{\eref}[1]	{(\protect\ref{#1})}

\providecommand{\paren}[1]	{\xmath{\left(#1\right)}}

\providecommand{\of}[1]		{\!\left(#1\right)}

\providecommand{\braces}[1]	{\xmath{\left\{#1\right\}}}

\providecommand{\xfrac}[2]	{\xmath{\frac{#1}{#2}}}

\providecommand{\Frac}[2]	{\xmath{{#1}/{#2}}}
\providecommand{\pfrac}[3][]	{\xmath{\!\paren{#1\frac{#2}{#3}}}}

\providecommand{\mathword}[3][]	{\,\inmath{\mathrm{#2}#1\of{#3}}}
\providecommand{\wordbrace}[3][]{\,\inmath{\mathrm{#2}#1\!\braces{#3}}}

\providecommand{\diag}[1]	{\wordbrace{diag}{#1}}

\providecommand{\sinc}		{\Ofun{\mathrm{sinc}}}

\providecommand{\sgn}[1] 	{\mathword{sgn}{#1}}

\usepackage{jf-accent}
\usepackage{jf-fun}
\usepackage{jf-names}

\newcommand{\st} {\xmath{\text{s.t.}\:}}

\usepackage{hyphenat}

\newcolumntype{B}{>{\centering\arraybackslash}p{4.15em}}
\newcolumntype{C}{>{\centering\arraybackslash}p{5.2em}}
\newcolumntype{A}{>{\centering\arraybackslash}p{1.6em}}
\newcolumntype{D}{>{\centering\arraybackslash}p{2.0em}}
\newcolumntype{E}{>{\centering\arraybackslash}p{5.4em}}

\newcommand{\cblue}[1] {\bgroup#1\egroup}
\renewcommand{\algref}[1] {Alg.~\ref{#1}\xspace}

\newcommand{\Reals} {\xmath{\mathbb{R}}}
\newcommand{\prox} {\operatorname{prox}} 
\newcommand{\FL} {\xmath{\mathcal{F}_{0,L}(\Reals^d)}}
\newcommand{\SL} {\xmath{\mathcal{F}_{\mu,L}(\Reals^d)}}
\newcommand{\Finf} {\xmath{\mathcal{F}_{0,\infty}(\Reals^d)}}
\newcommand{\betastar} {\xmath{\beta^{\star}}}
\newcommand{\gamstar} {\xmath{\gamma^{\star}}}
\newcommand{\bsig} {\xmath{\bar{\sigma}}}

\newcommand{\bxmath}[1] {\xmath{\bm{#1}}}
\newcommand{\x} {\bxmath{x}}
\newcommand{\y} {\bxmath{y}}
\newcommand{\z} {\bxmath{z}}
\newcommand{\w} {\bxmath{w}}
\newcommand{\vv} {\bxmath{v}}
\newcommand{\A} {\bxmath{A}}
\newcommand{\V} {\bxmath{V}}
\newcommand{\Lam} {\boldsymbol{\Lambda}}
\renewcommand{\a} {\bxmath{a}}
\newcommand{\bb} {\bxmath{b}}
\renewcommand{\lll} {\bxmath{l}}
\renewcommand{\u} {\bxmath{u}}
\newcommand{\Q} {\bxmath{Q}}
\newcommand{\vp} {\bxmath{p}}
\newcommand{\xii} {\boldsymbol{\xi}}
\newcommand{\I} {\bxmath{I}}
\newcommand{\Zero} {\bxmath{0}}
\newcommand{\T} {\bxmath{T}}
\newcommand{\omegaa} {\bxmath{\varepsilon}} 
\newcommand{\tr} {^\top} 

\usepackage{fullpage}
\usepackage{hyperref}
\hypersetup{
	linktoc=page,
	colorlinks=true,
	linkcolor=blue,
	citecolor=red,
	filecolor=magenta,
	urlcolor=cyan
}

\begin{document}

\epstopdfsetup{outdir=./} 

\title{
Adaptive Restart of the Optimized Gradient Method
for Convex Optimization
}


\author{Donghwan Kim         \and
        Jeffrey A. Fessler 
}


\institute{Donghwan Kim \and Jeffrey A. Fessler \at
                Dept. of Electrical Engineering and Computer Science,
                University of Michigan, Ann Arbor, MI 48109 USA \\
                \email{kimdongh@umich.edu, fessler@umich.edu}           
}

\date{Date of current version: \today} 

\maketitle

\begin{abstract}

First-order methods with momentum
such as Nesterov's fast gradient method
are very useful
for convex optimization problems,
but can exhibit undesirable oscillations
yielding slow convergence rates
for some applications.
An adaptive restarting scheme
can improve the convergence rate of the fast gradient method,
when the parameter
of a strongly convex cost function is unknown
or when the iterates of the algorithm
enter a locally strongly convex region.
Recently, we introduced the optimized gradient method,
a first-order algorithm
that has an inexpensive per-iteration computational cost
similar to that of the fast gradient method,
yet has a worst-case cost function rate
that is twice faster than that of the fast gradient method
and that is optimal for large-dimensional smooth convex problems.
Building upon the
success of accelerating the fast gradient method using adaptive restart,
this paper investigates similar heuristic acceleration of 
the optimized gradient method.
We first derive 
\cblue{a new first-order method 
that resembles the optimized gradient method}
for strongly convex quadratic problems
with known function parameters,
yielding a linear convergence rate
that is faster than that of the analogous version of 
the fast gradient method.
We then provide a heuristic analysis
and numerical experiments
that illustrate that adaptive restart
can accelerate the convergence of the optimized gradient method.
Numerical results also illustrate
that adaptive restart
is helpful
for a proximal version of the optimized gradient method
for nonsmooth composite convex functions.
\end{abstract}

\section{Introduction}

The computational expense of first-order methods
depends only mildly on the problem dimension,
so they are attractive
for solving large-dimensional optimization problems~\cite{cevher:14:cof}.
In particular, Nesterov's fast gradient method (FGM)
\cite{nesterov:83:amf,nesterov:04,beck:09:afi}
is used widely
because it has a worst-case cost function rate
that is optimal up to constant
for large-dimensional smooth convex problems
\cite{nesterov:04}.
In addition, for smooth and strongly convex problems
where the strong convexity parameter is known,
a version of FGM has a linear convergence rate~\cite{nesterov:04}
that improves upon that of a standard gradient method.
However,
without knowledge of the function parameters,
conventional
FGM does not guarantee a linear convergence rate.

When the strong convexity parameter is unknown,
a simple adaptive restarting scheme~\cite{odonoghue:15:arf}
for FGM
heuristically improves its convergence rate
(see also~\cite{giselsson:14:mar,su:16:ade} for theory
and~\cite{cevher:14:cof,muckley:15:fpm,monteiro:16:aaa} for applications).
In addition,
adaptive restart is useful
even when the function is only locally strongly convex
near the minimizer
\cite{odonoghue:15:arf}.
First-order methods are known to be suitable
when only moderate solution accuracy is required,
and
adaptive restart can help
first-order methods achieve medium to high accuracy.

Recently we proposed the optimized gradient method (OGM)~\cite{kim:16:ofo}
(built upon \cite{drori:14:pof})
that has efficient per-iteration computation similar to FGM
yet that exactly achieves the optimal worst-case rate
for decreasing a large-dimensional smooth convex function
among all first-order methods
\cite{drori:17:tei}.
(See~\cite{kim:16:gto-arxiv,kim:17:otc,taylor:17:ewc}
for further analysis and extensions of OGM.)
This paper examines 
\cblue{a general class of accelerated first-order methods 
that includes a gradient method (GM), FGM, and OGM}
for strongly convex \emph{quadratic} functions,
and 
develops an OGM variant, \cblue{named OGM-$q$},
that provides a linear convergence rate
that is faster than that of \cblue{the analogous version of} FGM.
The analysis reveals that,
like FGM~\cite{odonoghue:15:arf},
OGM may exhibit undesirable oscillating behavior
in some cases.
Building on the quadratic analysis 
and the adaptive restart scheme of FGM in~\cite{odonoghue:15:arf},
we propose an adaptive restart scheme 
that heuristically accelerates the convergence rate of OGM
when the function is strongly convex 
or even when it is only locally strongly convex.
This restart scheme
circumvents the oscillating behavior.
Numerical results illustrate
that the proposed OGM with restart
performs better than
FGM with restart in~\cite{odonoghue:15:arf}.

Sec.~\ref{sec:prob,algo}
reviews first-order methods
for convex problems
such as GM, FGM, and OGM.
Sec.~\ref{sec:quadanal}
\cblue{analyzes a general class of accelerated first-order methods
that includes GM, FGM, and OGM}
for strongly convex quadratic problems,
\cblue{and proposes a new OGM variant
with a fast linear convergence rate}.
Sec.~\ref{sec:restart}
suggests an adaptive restart scheme for OGM
using the quadratic analysis in Sec.~\ref{sec:quadanal}.
Sec.~\ref{sec:propOGM}
illustrates the proposed adaptive version of OGM
that we use for numerical experiments 
on various convex problems in Sec.~\ref{sec:result},
including nonsmooth composite convex functions, 
and
Sec.~\ref{sec:conc} concludes.

\section{Problem and Methods}
\label{sec:prob,algo}

\subsection{Smooth and Strongly Convex Problem}

We first consider
the smooth and strongly convex minimization problem:
\begin{align}
\min_{\x\in\Reals^d} \;&\; f(\x)
\label{eq:prob} \tag{M}
\end{align}
that satisfies the following smooth and strongly convex conditions:

\begingroup
\allowdisplaybreaks
\begin{itemize} 
\item
$f\;:\;\Reals^d\rightarrow\Reals$
has Lipschitz continuous gradient with Lipschitz constant $L>0$, \ie,
\begin{align}
||\nabla f(\x) - \nabla f(\y)|| \le L||\x-\y||, \quad \forall \x, \y\in\Reals^d
,\end{align}
\item $f$ is strongly convex with strong convexity parameter $\mu>0$, \ie,
\begin{align}
f(\x) \ge f(\y) + \Inprod{\nabla f(\y)}{\x - \y} + \frac{\mu}{2}||\x - \y||^2,
 \quad \forall \x, \y\in\Reals^d
\label{eq:strcvx}
.\end{align}
\end{itemize}
\endgroup
We let \SL denote the class of functions $f$
that satisfy the above two conditions hereafter,
and let $\x_*$ denote the unique minimizer of $f$.
We let $q := \Frac{\mu}{L}$
denote the reciprocal of the condition number
of a function $f \in \SL$.
We also let \FL 
denote the class of smooth convex functions $f$
that satisfy the above two conditions with $\mu=0$,
and let $\x_*$ denote a minimizer of $f$.

Some algorithms discussed in this paper
require knowledge of both $\mu$ and $L$,
but in many cases estimating $\mu$ is challenging compared to computing $L$.\footnote{
For some applications even estimating $L$ is expensive,
and one must employ a backtracking scheme~\cite{beck:09:afi}
or similar approaches.
We assume $L$ is known throughout this paper.
An estimate of $\mu$ could be found
by a backtracking scheme
as described in~\cite[Sec. 5.3]{nesterov:13:gmf}.
}
Therefore,
this paper focuses on the case where
the parameter $\mu$ is unavailable while $L$ is available.
Even without knowing $\mu$,
the adaptive restart approach in~\cite{odonoghue:15:arf}
and the proposed adaptive restart approach in this paper
both exhibit linear convergence rates
in strongly convex cases.

We next review known accelerated first-order methods
for solving~\eqref{eq:prob}.

\subsection{\cblue{Review of} Accelerated First-order Methods}

This paper focuses on accelerated first-order \cblue{methods (AFM)} 
of the form shown in~\algref{alg:ogm}.
The fast gradient method (FGM)
\cite{nesterov:83:amf,nesterov:04,beck:09:afi}
(with $\gamma_k = 0$ in \algref{alg:ogm})
accelerates the gradient method (GM) 
(with $\beta_k = \gamma_k = 0$)
using the \emph{momentum} term $\beta_k(\y_{k+1} - \y_k)$
with negligible additional computation.
The optimized gradient method (OGM)
\cite{kim:16:ofo,kim:17:otc}
uses
an \emph{over-relaxation} term
$\gamma_k(\y_{k+1} - \x_k) = -\gamma_k\alpha\nabla f(\x_k)$
for further acceleration.

\begin{algorithm}[H]
\caption{Accelerated First-order \cblue{Methods (AFM)}}
\label{alg:ogm}
\begin{algorithmic}[1]
\State {\bf Input:} $f\in\FL$ or $\SL$, $\x_0 = \y_0\in\Reals^d$.
\For{$k \ge 0$}
\State $\y_{k+1} = \x_k - \alpha\nabla f(\x_k)$
\State $\x_{k+1} = \y_{k+1}
        + \beta_k(\y_{k+1} - \y_k)
        + \gamma_k(\y_{k+1} - \x_k)$
\EndFor
\end{algorithmic}
\end{algorithm}

\cblue{
Tables~\ref{tab:alg1} and~\ref{tab:alg2}
summarize
the standard choices of coefficients $(\alpha,\beta_k,\allowbreak \gamma_k)$
for GM, FGM, OGM
in~\cite{nesterov:83:amf,nesterov:04,beck:09:afi,kim:16:ofo,kim:17:otc}
and their worst-case rates
for smooth convex functions $\FL$
and smooth and strongly convex functions $\SL$
respectively.
(Other choices can be found in
\cite{nesterov:04,kim:16:gto-arxiv,chambolle:15:otc}.)
For convenience hereafter,
we use the names GM, GM-$q$, FGM, FGM-$q$, OGM, and OGM$'$
to distinguish different choices of standard AFM coefficients
in Tables~\ref{tab:alg1} and~\ref{tab:alg2}.
}

\begin{table}[htbp] 
\centering
\renewcommand{\arraystretch}{1.4} 
\caption{ 
\cblue{Accelerated First-order Methods for Smooth Convex Problems}}
\label{tab:alg1}
\cblue{
\begin{tabular}{|B||A|A|A||c|}
\hline
Method & $\alpha$ & $\beta_k$ & $\gam_k$ & Worst-case Rate \\
\hline
GM  & $\frac{1}{L}$
        & $0$
        & $0$
        & $\scriptstyle f(\y_k) - f(\x_*) \le \frac{L||\x_0-\x_*||^2}{4k+2}$
	\cite{drori:14:pof}
        \\ \hline
\multirow{2}{*}{FGM~\cite{nesterov:83:amf}} 
	& \multirow{2}{*}{$\frac{1}{L}$}
        & \multirow{2}{*}{$\frac{t_k-1}{t_{k+1}}$}
        & \multirow{2}{*}{$0$}
        & $\scriptstyle f(\y_k) - f(\x_*) \le \frac{L||\x_0-\x_*||^2}{2t_{k-1}^2}
			\le \frac{2L||\x_0-\x_*||^2}{(k+1)^2}$ 
	\cite{beck:09:afi} \\
	& & &
	& $\scriptstyle f(\x_k) - f(\x_*) \le \frac{L||\x_0-\x_*||^2}{2t_k^2}
			\le \frac{2L||\x_0-\x_*||^2}{(k+2)^2}$ 
	\cite{kim:16:ofo}
        \\ \hline
OGM$'$~\cite{kim:17:otc}
	& $\frac{1}{L}$
        & $\frac{t_k-1}{t_{k+1}}$
        & $\frac{t_k}{t_{k+1}}$
        & $\scriptstyle f(\y_k) - f(\x_*) \le \frac{L||\x_0-\x_*||^2}{4t_{k-1}^2}
			\le \frac{L||\x_0-\x_*||^2}{(k+1)^2}$
        \cite{kim:17:otc}
	\\ \hline
OGM~\cite{kim:16:ofo}
        & $\frac{1}{L}$
        & $\frac{\theta_k-1}{\theta_{k+1}}$
        & $\frac{\theta_k}{\theta_{k+1}}$
        & $\scriptstyle f(\x_N) - f(\x_*) \le \frac{L||\x_0-\x_*||^2}{2\theta_N^2}
			\le \frac{L||\x_0 - \x_*||^2}{(N+1)^2}$
        \cite{kim:16:ofo}
	\\ \hline \hline
\multicolumn{5}{|c|}
{Parameters}
	\\ \hline
\multicolumn{5}{|c|}
{
\hspace{-4.4em}
$\scriptstyle
\hspace{0.8em}
t_0 = 1, \;\;\;
t_k = \frac{1}{2}\paren{1+\sqrt{1+4t_{k-1}^2}},\;\; k=1,\ldots,$ 
} \\
\multicolumn{5}{|c|}
{
$\scriptstyle
\theta_0 = 1, \;\;
\theta_k = \begin{cases}
	\scriptstyle
	\frac{1}{2}\paren{1+\sqrt{1+4\theta_{k-1}^2}}, & \scriptstyle k=1,\ldots,N-1, \\
	\scriptstyle
	\frac{1}{2}\paren{1+\sqrt{1+8\theta_{k-1}^2}}, & \scriptstyle k=N.
	\end{cases}
$
} \\ \hline
\end{tabular}
}
\end{table}
\vspace{-5pt}

\begin{table}[htbp] 
\centering
\renewcommand{\arraystretch}{1.4} 
\caption{ 
\cblue{
Accelerated First-order Methods (with $\gamma_k = 0$) 
for Smooth and Strongly Convex Problems}
\cblue{(The worst-case rates also apply to $\frac{\mu}{2}||\y_k - \x_*||^2$
due to the strong convexity~\eqref{eq:strcvx}.)}}
\label{tab:alg2}
\cblue{
\begin{tabular}{|C||D|D||c|}
\hline
Method & $\alpha$ & $\beta_k$ & Worst-case Rate \\
\hline
GM 
	& $\frac{1}{L}$
        & $0$
        & $\scriptstyle f(\y_k) - f(\x_*) \le \paren{1 - \frac{2\mu}{1+q}}^k\frac{L||\x_0-\x_*||^2}{2}$
	\cite{nesterov:04}
        \\ \hline 
GM-$q$
	& $\frac{2}{\mu+L}$
        & $0$
        & $\scriptstyle f(\y_k) - f(\x_*) \le \paren{\frac{1-q}{1+q}}^{2k}\frac{L||\x_0-\x_*||^2}{2}$
	\cite{nesterov:04}
        \\ \hline
FGM-$q$ \cite{nesterov:04}
	& $\frac{1}{L}$
        & $\frac{1-\sqrt{q}}{1+\sqrt{q}}$
        & $\scriptstyle f(\y_k) - f(\x_*) \le (1-\sqrt{q})^k\frac{(1+q)L||\x_0-\x_*||^2}{2}$
        \cite{nesterov:04}
	\\ \hline
\end{tabular}
}
\end{table}
\vspace{-5pt}

\cblue{
The worst-case OGM rate~\cite{kim:16:ofo} in Table~\ref{tab:alg1} 
is about twice faster than the FGM rate~\cite{beck:09:afi}}
and is optimal for first-order methods
for the function class \FL
under the large-scale condition
$d\ge N+1$~\cite{drori:17:tei}.
\cblue{However,}
it is yet unknown
\cblue{which first-order methods provide an optimal worst-case} 
linear convergence rate 
for the function class \SL;
this topic is left as an interesting future work.\footnote{
Recently,
\cite{vanscoy:18:tfk}
developed a new first-order method for known $q$
that \cblue{is not in AFM class but} 
achieves a linear worst-case rate $(1-\sqrt{q})^2$ 
for the decrease of a strongly convex function
that is faster than the linear rate $(1-\sqrt{q})$
\cblue{of FGM-$q$ in Table~\ref{tab:alg2}.}}
Towards this direction,
Sec.~\ref{sec:quadanal}
studies \cblue{AFM} for strongly convex \emph{quadratic} problems,
\cblue{leading to a new method named OGM-$q$
with a linear convergence rate that is faster than that of FGM-$q$}.
Sec.~\ref{sec:restart} uses
this quadratic analysis
to analyze an adaptive restart scheme for OGM.

\section{
Analysis of \cblue{AFM} for Quadratic Functions
}
\label{sec:quadanal}

This section analyzes the behavior of \cblue{AFM}
for minimizing a strongly convex quadratic function.
The quadratic analysis of \cblue{AFM} 
in this section is similar
in spirit to the analyses of a heavy-ball method~\cite[Sec. 3.2]{polyak:87}
and \cblue{AFM with $\gamma_k=0$}
\cite[Appx. A]{lessard:16:aad}
\cite[Sec. 4]{odonoghue:15:arf}.

\cblue{In addition, Sec.~\ref{sec:ogmq}}
optimizes the coefficients of \cblue{AFM} for such quadratic functions,
yielding a linear convergence rate that is faster than that of
\cblue{FGM-$q$}.
The resulting \cblue{method, named OGM-$q$,} requires the knowledge of $q$,
and \cblue{Sec.~\ref{sec:OGM,conv}} 
shows that using 
\cblue{OGM (and OGM$'$) in Table~\ref{tab:alg1}}
instead
(without the knowledge of $q$)
will cause the OGM iterates to oscillate
when the momentum is larger than a critical value.
This analysis stems
from the dynamical system analysis of 
\cblue{AFM with $\alpha=\Frac{1}{L}$ and $\gamma_k=0$} 
in~\cite[Sec. 4]{odonoghue:15:arf}.

\subsection{Quadratic Analysis of \cblue{AFM}}
\label{sec:quad}

This section
considers minimizing a strongly convex quadratic function:
\begin{align}
f(\x) = \frac{1}{2}\x^\top \Q \x - \vp\tr \x \in \SL
\label{eq:quad}
\end{align}
where $\Q \in \Reals^{d\times d}$
is a symmetric positive definite matrix,
$\vp \in \Reals^d$ is a vector.
Here, $\nabla f(\x) = \Q\x - \vp$ is the gradient,
and $\x_* = \Q^{-1} \vp$ is the optimum.
The smallest and the largest eigenvalues of \Q
correspond to the parameters $\mu$ and $L$ of the function respectively.
For simplicity in the quadratic analysis,
we consider the version of \cblue{AFM}
that has constant coefficients
$(\alpha,\beta,\gamma)$.

Defining the vectors
$\xii_k := (\x_k^\top, \x_{k-1}^\top)^\top \in \Reals^{2d}$
and
$\xii_* := (\x_*^\top, \x_*^\top)^\top\in \Reals^{2d}$,
and extending the analysis for \cblue{AFM with $\gamma=0$}
in~\cite[Appx. A]{lessard:16:aad},
\cblue{AFM} 
has
the following equivalent form
for $k \ge 1$:
\begin{align}
&\xii_{k+1} - \xii_* = \T(\alpha,\beta,\gamma) \, (\xii_k - \xii_*)
\label{eq:xii,update}
,\end{align}
where the system matrix
$\T(\alpha,\beta,\gamma)$ of \cblue{AFM} is defined as
\begin{align}
\T(\alpha,\beta,\gamma) := \bigg[\begin{array}{cc}
	(1+\beta)(\I-\alpha\Q) - \gamma\alpha\Q & -\beta(\I-\alpha\Q) \\
	\I & \Zero
	\end{array}\bigg] \in \Reals^{2d\times 2d}
\label{eq:T}
\end{align}
for an identity matrix $\I \in \Reals^{d\times d}$.
The sequence
$\{\tilde{\xii}_k := (\y_k^\top, \y_{k-1}^\top)^\top\}_{k\ge1}$
also satisfies the recursion~\eqref{eq:xii,update},
implying that~\eqref{eq:xii,update}
characterizes the behavior of both the primary sequence $\{\y_k\}$
and the secondary sequence $\{\x_k\}$ of \cblue{AFM}
with constant coefficients.

The spectral radius $\rho(\T(\cdot))$ of matrix $\T(\cdot)$
determines the convergence rate of the algorithm.
Specifically,
for any $\epsilon>0$,
there exists $K\ge0$ such that
$[\rho(\T)]^k \le ||\T^k|| \le (\rho(\T) + \epsilon)^k$
for all $k\ge K$,
establishing
the following worst-case rate:
\begin{align}
||\xii_{k+1} - \xii_*||^{\cblue{2}} \le
(\rho(\T(\alpha,\beta,\gamma)) + \epsilon)^{\cblue{2}k} 
\ ||\xii_1 - \xii_*||^{\cblue{2}}
\label{eq:xii,rate}
.\end{align}
We next analyze
$\rho(\T(\alpha,\beta,\gamma))$.

Considering the eigen-decomposition of \Q in $\T(\cdot)$
as in~\cite[Appx. A]{lessard:16:aad},
the spectral radius of $\T(\cdot)$
is:
\begin{align}
\rho(\T(\alpha,\beta,\gamma))
	= \max_{\mu\le\lambda\le L} \rho(\T_\lambda(\alpha,\beta,\gamma))
\label{eq:rhoT}
,\end{align}
where
for any eigenvalue $\lambda$ 
of matrix \Q we define 
a matrix $\T_\lambda(\alpha,\beta,\gamma) \in \Reals^{2\times 2}$
by substituting $\lambda$ and $1$ for $\Q$ and $\I$ in 
$\T(\alpha,\beta,\gamma)$ respectively.
Similar to the analysis of \cblue{AFM with $\gamma=0$} 
in~\cite[Appx. A]{lessard:16:aad},
the spectral radius of $\T_\lambda(\alpha,\beta,\gamma)$
is:
\begingroup
\allowdisplaybreaks
\begin{align}
&\, \rho(\T_\lambda(\alpha,\beta,\gamma))
= \max\{|r_1(\alpha,\beta,\gamma,\lambda)|, |r_2(\alpha,\beta,\gamma,\lambda)|\}
	\label{eq:rho} \\
=&\, \begin{cases}
	\frac{1}{2} \paren{
	|(1 + \beta)\paren{1 - \alpha\lambda} - \gamma \alpha\lambda|
	+ \sqrt{\Delta(\alpha,\beta,\gamma,\lambda)} },
	& \Delta(\alpha,\beta,\gamma,\lambda) \ge 0, \\
	\sqrt{\beta(1-\alpha\lambda)}, & \text{otherwise,}
	\nonumber
\end{cases}
\end{align}
\endgroup
where $r_1(\alpha,\beta,\gamma,\lambda)$
and $r_2(\alpha,\beta,\gamma,\lambda)$
denote the roots of the characteristic polynomial of $\T_\lambda(\cdot)$:
\begin{align}
r^2 - ((1+\beta)(1-\alpha\lambda) - \gamma\alpha\lambda)r + \beta(1-\alpha\lambda)
\label{eq:poly}
,\end{align}
and
$\Delta(\alpha,\beta,\gamma,\lambda)
	:= \paren{(1 + \beta)\paren{1 - \alpha\lambda} - \gamma \alpha\lambda}^2
	- 4 \beta\paren{1 - \alpha\lambda}$
denotes the corresponding discriminant.
For fixed $(\alpha,\beta,\gamma)$,
the spectral radius
$\rho(\T_\lambda(\alpha,\beta,\gamma))$
in~\eqref{eq:rho}
is a continuous and quasi-convex\footnote{
\label{ftquasi}
It is straightforward to show that
$\rho(\T_\lambda(\alpha,\beta,\gamma))$ in~\eqref{eq:rho}
is quasi-convex over $\lambda$.
First, $\sqrt{\beta(1-\alpha\lambda)}$ is quasi-convex over $\lambda$
(for $\Delta(\alpha,\beta,\gamma,\lambda) < 0$).
Second, the eigenvalue $\lambda$ satisfying
$\Delta(\alpha,\beta,\gamma,\lambda) \ge 0$
is in the region where
the function
$ 
\frac{1}{2} \paren{
        |(1 + \beta)\paren{1 - \alpha\lambda} - \gamma \alpha\lambda|
        + \sqrt{\Delta(\alpha,\beta,\gamma,\lambda)} }
$
either monotonically increases or decreases,
which overall makes the continuous function $\rho(\T_\lambda(\alpha,\beta,\gamma))$
quasi-convex over $\lambda$.
This proof can be simply applied to other variables, \ie,
$\rho(\T_\lambda(\alpha,\beta,\gamma))$
is quasi-convex over either $\alpha$, $\beta$ or \gam.
} 
function of $\lambda$;
thus its maximum
over $\lambda$
occurs at one of its boundary points
$\lambda = \mu$ or $\lambda = L$.


\cblue{The next section reviews the optimization of AFM coefficients
to provide the fastest convergence rate,
\ie, the smallest spectral radius $\rho(\T(\cdot))$ in \eqref{eq:rhoT},
under certain constraints on $(\alpha,\beta,\gamma)$.
}

\subsection{\cblue{Review of Optimizing AFM} Coefficients
\cblue{under Certain Constraints on $(\alpha,\beta,\gamma)$}}

\cblue{The AFM} coefficients
that provide the fastest convergence
for minimizing a strongly convex quadratic function
\cblue{would} solve
\begin{align}
\argmin{\alpha,\beta,\gamma} \rho(\T(\alpha,\beta,\gamma))
=
\argmin{\alpha,\beta,\gamma}
	\max\{\rho(\T_{\mu}(\alpha,\beta,\gamma)),\rho(\T_L(\alpha,\beta,\gamma))\}
\label{eq:optogm}
.\end{align}
\cblue{Note that a heavy-ball method~\cite{polyak:87}
(that is not in AFM class)
with similarly optimized coefficients
has a linear worst-case rate 
with $\rho(\cdot) = \frac{1 - \sqrt{q}}{1 + \sqrt{q}}$
that is optimal (up to constant)
for strongly convex quadratic problems
\cite{nesterov:04}.
Thus, optimizing~\eqref{eq:optogm} 
would be of little practical benefit for quadratic problems.
Nevertheless, such optimization is new to AFM for $\gamma>0$ 
(with the additional constraint $\alpha=\Frac{1}{L}$ introduced below),
and is useful in our later analysis for the adaptive restart
in Sec.~\ref{sec:restart}.
A heavy-ball method with the coefficients 
optimized for strongly convex quadratic problems
does not converge
for some strongly convex nonquadratic problems~\cite{lessard:16:aad},
and other choices of coefficients
do not yield worst-case rates that are comparable 
to those of some accelerated choices of AFM
\cite{drori:14:pof,lessard:16:aad}, 
so we focus on AFM hereafter.}

\cblue{The coefficient optimization~\eqref{eq:optogm} for AFM} 
was studied previously
\cblue{with various constraint}.
\cblue{For example,} 
optimizing~\eqref{eq:optogm} over $\alpha$ 
\cblue{with the constraint $\beta=\gamma=0$
yields GM-$q$.}
Similarly,
\cblue{FGM-$q$}
results from
optimizing~\eqref{eq:optogm} over $\beta$
for the \cblue{constraint}%
\footnote{%
For \cblue{FGM-$q$} 
the value of 
$\rho(\T_L(\Frac{1}{L},\beta,0))$ is $0$,
and the function $\rho(\T_{\mu}(\Frac{1}{L},\beta,0))$
is continuous and quasi-convex over $\beta$
(see footnote~\ref{ftquasi}).
The minimum of $\rho(\T_{\mu}(\Frac{1}{L},\beta,0))$
occurs at the point $\beta = \frac{1-\sqrt{q}}{1+\sqrt{q}}$
in \cblue{Table~\ref{tab:alg2}}
satisfying $\Delta\paren{\Frac{1}{L},\beta,0,\mu} = 0$,
verifying the statement
that \cblue{FGM-$q$} 
results from optimizing~\eqref{eq:optogm} over $\beta$
given $\alpha = \Frac{1}{L}$ and $\gamma = 0$.
}
$\alpha = \Frac{1}{L}$ and $\gamma = 0$.
\cblue{In~\cite[Prop.~1]{lessard:16:aad},
AFM with coefficients
\(
(\alpha,\beta,\gamma) 
= \paren{\frac{4}{\mu + 3L},
\frac{\sqrt{3+q} - 2\sqrt{q}}{\sqrt{3+q} + 2\sqrt{q}},\,0}
,\)
named FGM$'$-$q$ in Table~\ref{tab:coeff},
was derived by optimizing~\eqref{eq:optogm} over $(\alpha,\beta)$
with the constraint $\gamma=0$.}

Although a general unconstrained solution to~\eqref{eq:optogm}
would be an interesting future direction,
here
we focus on optimizing~\eqref{eq:optogm} over $(\beta,\gamma)$
\cblue{with the constraint} $\alpha=\Frac{1}{L}$.
This choice
simplifies the problem~\eqref{eq:optogm}
and is useful for analyzing
an adaptive restart scheme for OGM in Sec.~\ref{sec:restart}.

\subsection{Optimizing the Coefficients $(\beta,\gamma)$ of \cblue{AFM}
When $\alpha = \Frac{1}{L}$}
\label{sec:ogmq}

When $\alpha = \Frac{1}{L}$ and $\lambda = L$,
the characteristic polynomial~\eqref{eq:poly}
becomes
\(
r^2 + \gamma r = 0
.\)
The roots are $r = 0$ and $r = -\gamma$,
so $\rho(\T_L(\Frac{1}{L},\beta,\gamma)) = |\gamma|$.
In addition, because
$\rho(\T_{\mu}(\Frac{1}{L},\beta,\gamma))$
is continuous and quasi-convex over $\beta$
(see footnote~\ref{ftquasi}),
it can be easily shown that
the smaller value of $\beta$
satisfying the following equation:
\begin{align}
&\,\Delta(\Frac{1}{L},\beta,\gamma,\mu)
= ((1+\beta)(1-q) - \gamma q)^2 - 4\beta(1-q) \\
=&\, (1-q)^2\beta^2 - 2(1-q)(1+q+q\gamma)\beta + (1-q)(1-q-2q\gamma) + q^2\gamma^2
	= 0 \nonumber
\end{align}
minimizes $\rho(\T_{\mu}(\Frac{1}{L},\beta,\gamma))$
for any given \gam (satisfying \cblue{$\gamma \ge -1$}). 
The optimal
$\beta$ \cblue{for a given $\gamma$ (when $\alpha=\Frac{1}{L}$)}
is
\begin{align}
&\betastar(\gamma) := \Frac{\paren{1 - \sqrt{q(1 + \gamma)}}^2}{(1-q)}
\label{eq:betaopt}
,\end{align}
which reduces to
$\beta = \betastar(0) = \frac{1-\sqrt{q}}{1+\sqrt{q}}$
for \cblue{FGM-$q$} (with $\gamma = 0$).
Substituting
\eqref{eq:betaopt} into~\eqref{eq:rho}
yields
\(
\rho(\T_{\mu}(\Frac{1}{L},\betastar(\gamma),\gamma))
= |1 - \sqrt{q(1+\gamma)}|
,\)
leading
to the following
simplification of~\eqref{eq:optogm}
with $\alpha = \Frac{1}{L}$ and $\beta = \betastar(\gamma)$
from~\eqref{eq:betaopt}:
\begin{align}
\gamstar := \argmin{\gamma} \,
\max\left\{ |1 - \sqrt{q(1 + \gamma)} |, \, |\gamma| \right\}
\label{eq:optogmalpbet}
.\end{align}
The minimizer of~\eqref{eq:optogmalpbet}
satisfies
\(
1 - \sqrt{q(1 + \gamma)} = \pm \gamma
,\)
and with simple algebra,
we get
the following solutions to~\eqref{eq:optogm} 
with \cblue{the constraint} $\alpha = \Frac{1}{L}$ 
(and~\eqref{eq:optogmalpbet}):
\begin{align}
\betastar := \beta^\star(\gamstar)
	= \frac{\paren{\gamstar}^2}{1-q}
	= \frac{(2+q - \sqrt{q^2 + 8q})^2}{4(1-q)},
\quad
\gamstar = \frac{2 + q - \sqrt{q^2 + 8q}}{2}
\label{eq:ogmquadalp}
,\end{align}
for which
the spectral radius is
\(
\rho^\star
:= \rho( \T(\Frac{1}{L}, \betastar, \gamstar) )
= 1 - \sqrt{q(1+\gamstar)}
= \gamstar
.\)
\cblue{We denote \algref{alg:ogm} with coefficients 
$\alpha=\Frac{1}{L}$ and (\betastar,\gamstar) in~\eqref{eq:ogmquadalp}
as OGM-$q$.}

Table~\ref{tab:coeff}
compares
the spectral radius
of the 
\cblue{OGM-$q$}
to \cblue{GM-$q$, FGM-$q$, and FGM$'$-$q$~\cite[Prop.~1]{lessard:16:aad}}.
Simple algebra
shows that the spectral radius of 
\cblue{OGM-$q$} is smaller than those of \cblue{FGM-$q$ and FGM$'$-$q$},
\ie,
\(
\frac{2 + q - \sqrt{q^2 + 8q}}{2}
 \le 1 - \frac{2\sqrt{q}}{\sqrt{3+q}}
 \le 1 - \sqrt{q}
.\)
Therefore,
\cblue{OGM-$q$} 
achieves a worst-case convergence rate
of $||\xii_k - \xii_*||^{\cblue{2}}$
that is faster than that of FGM \cblue{variants}
\cblue{(but that is slower than a heavy-ball method~\cite{polyak:87})}
for a strongly convex quadratic function.

\newcommand{\mycol}[1] { \colorbox{gray!25}{ \makebox[3em]{ $ #1 $ }}}

\begin{table}[htbp] 
\centering
\renewcommand{\arraystretch}{1.4} 
\caption{
Optimally tuned coefficients $(\alpha,\beta,\gamma)$
of \cblue{GM-$q$, FGM-$q$, FGM$'$-$q$, and OGM-$q$},
and their spectral radius $\rho(\T(\alpha,\beta,\gamma))$~\eqref{eq:rhoT}.
These optimal coefficients result from solving~\eqref{eq:optogm}
with the shaded coefficients fixed.}
\label{tab:coeff}
\begin{tabular}{|C||E|E|E||E|}
\hline
Method & $\alpha$ & $\beta$ & \gam & $\rho(\T(\alpha,\beta,\gamma))$ \\
\hline
\vspace{-1.1em} \cblue{GM-$q$}  
	& $\frac{2}{\mu+L}$
	& \mycol{ 0 }
	& \mycol{ 0 }
	& $\frac{1-q}{1+q}$
	\\ \hline
\vspace{-1.1em} \cblue{FGM-$q$ \cite{nesterov:04}}
	& \mycol{ \frac{1}{L} }
	& $\frac{1-\sqrt{q}}{1+\sqrt{q}}$
	& \mycol{ 0 } 
	& $1-\sqrt{q}$
	\\ \hline 
\vspace{-1.1em} \cblue{FGM$'$-$q$ \cite{lessard:16:aad}}
	& $\frac{4}{\mu+3L}$ 
	& $\frac{\sqrt{3+q}-2\sqrt{q}}{\sqrt{3+q}+2\sqrt{q}}$
	& \mycol{ 0 }
	& $1 - \frac{2\sqrt{q}}{\sqrt{3+q}}$
	\\ \hline
\vspace{-1.1em} \cblue{OGM-$q$}
	& \mycol{ \frac{1}{L} }
	& $\frac{(2+q - \sqrt{q^2 + 8q})^2}{4(1-q)}$
	& $\frac{2+q - \sqrt{q^2 + 8q}}{2}$
	& $\frac{2+q - \sqrt{q^2 + 8q}}{2}$
	\\ \hline
\end{tabular}
\end{table}

\begin{figure}[htbp] 
\begin{center}
\includegraphics[clip,width=0.48\textwidth]{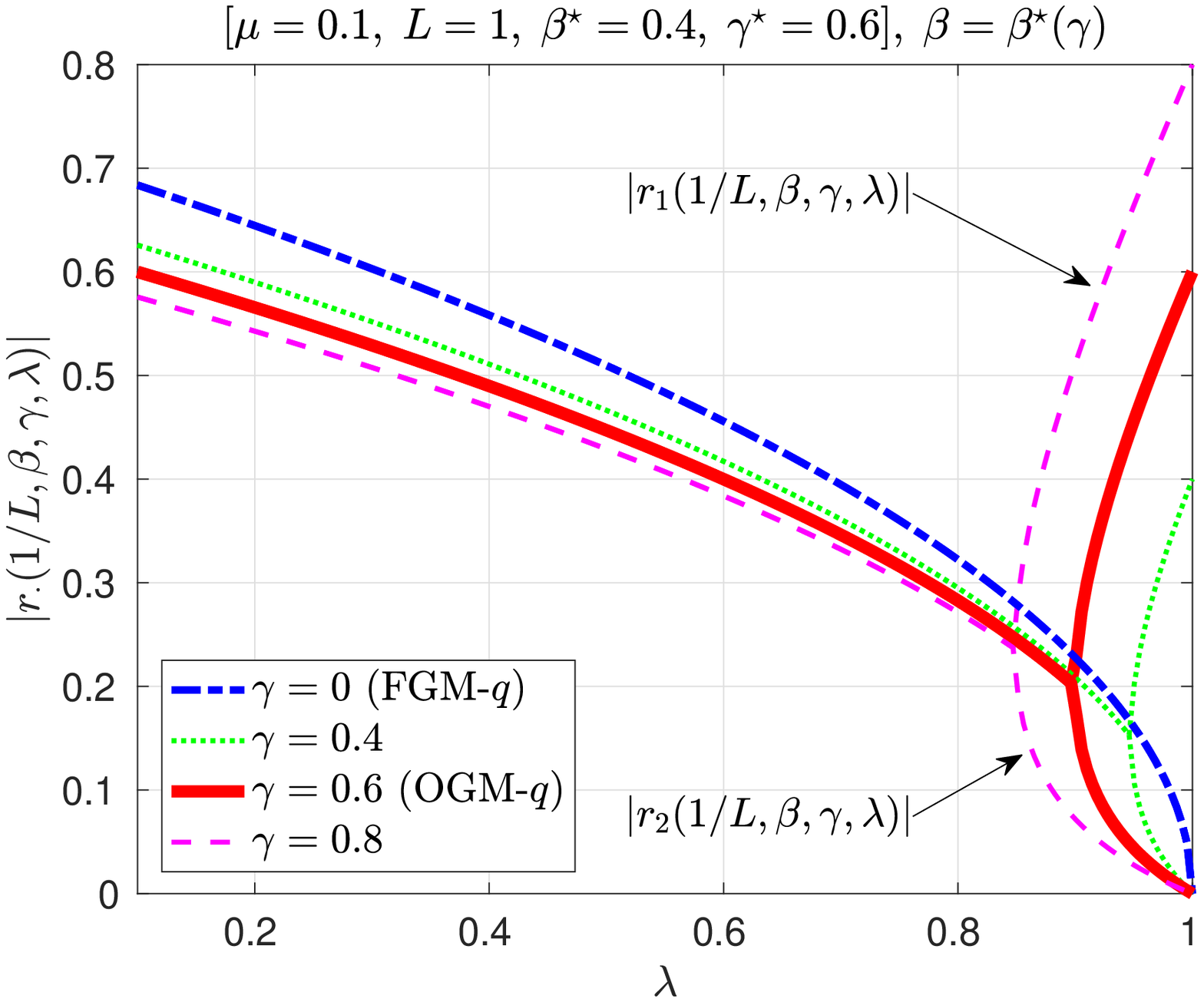}
\includegraphics[clip,width=0.48\textwidth]{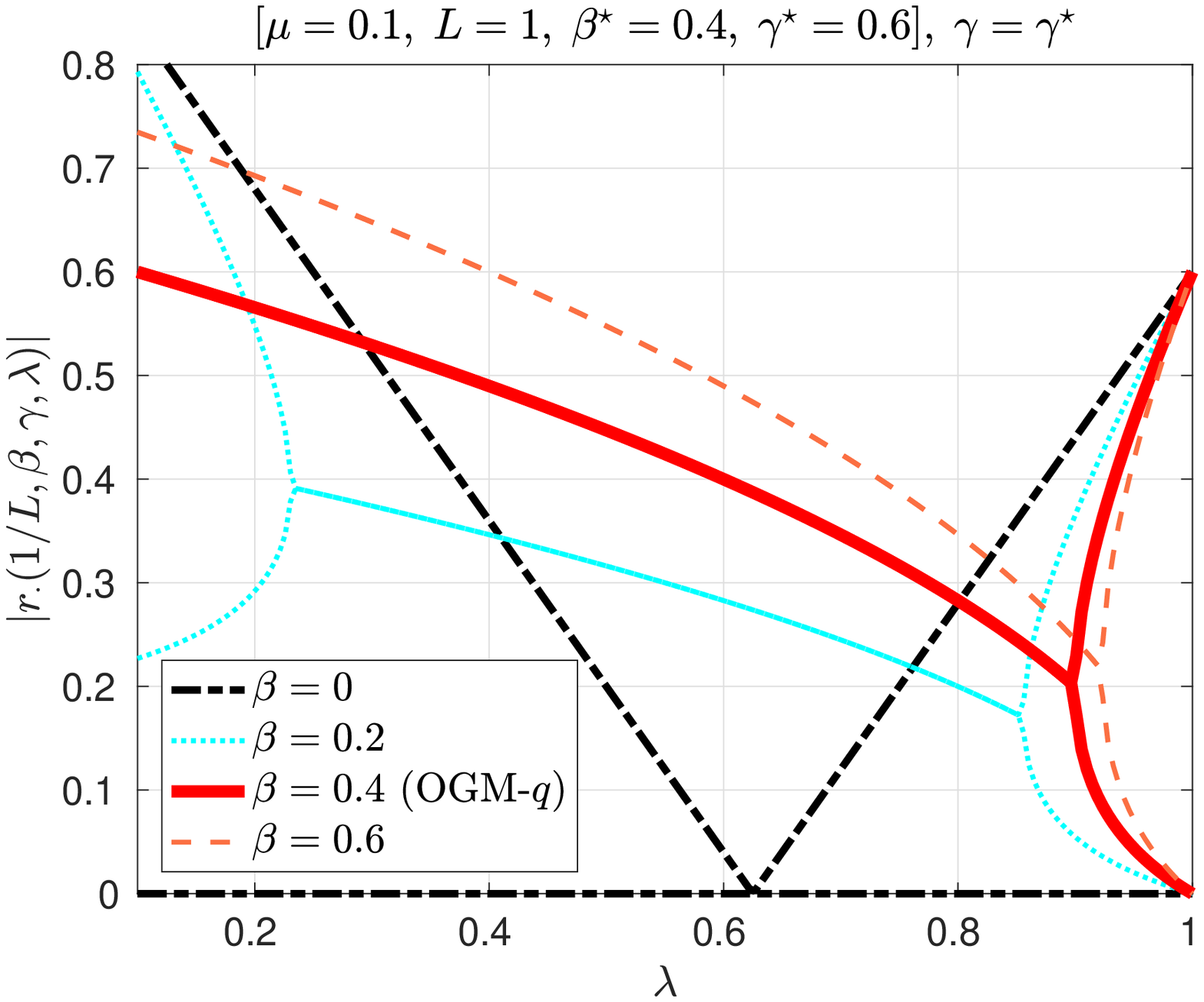}
\end{center}
\caption{
Plots of
$|r_1(\Frac{1}{L},\beta,\gamma,\lambda)|$
and
$|r_2(\Frac{1}{L},\beta,\gamma,\lambda)|$
over $\mu\le\lambda\le L$
for various
(Left) \gam values for given $\beta = \betastar(\gamma)$,
and
(Right) $\beta$ values for given $\gamma = \gamstar$,
for a strongly convex quadratic problem with 
$\mu=0.1$ and $L=1$ ($q=0.1$),
where $(\betastar,\gamstar) = (0.4,0.6)$.
%
The maximum of
$|r_1(\Frac{1}{L},\beta,\gamma,\lambda)|$
and
$|r_2(\Frac{1}{L},\beta,\gamma,\lambda)|$,
\ie the upper curve in the plot,
corresponds to 
the value of $\rho(\T_\lambda(\Frac{1}{L},\beta,\gamma))$ in~\eqref{eq:rho},
\cblue{and
the maximum value of $\rho(\T_\lambda(\Frac{1}{L},\beta,\gamma))$ over $\lambda$
corresponds to 
a spectral radius $\rho(\T(\Frac{1}{L},\beta.\gamma))$ in~\eqref{eq:rhoT}.}
}
\label{fig:ogm,roots}
\vspace{-3pt}
\end{figure}

To further understand the behavior of \cblue{AFM}
for each eigen-mode,
Fig.~\ref{fig:ogm,roots}
plots
$\rho(\T_\lambda(\Frac{1}{L},\beta,\gamma))$
\cblue{over} $\mu \le \lambda \le L$
for $\mu = 0.1$ and $L = 1$ ($q=0.1$) as an example,
where $(\betastar,\gamstar) = (0.4,0.6)$.
\cblue{The left plot of} Fig.~\ref{fig:ogm,roots}
first compares 
\cblue{the $\rho(\T_{\lambda}(\Frac{1}{L},\beta,\gamma))$ values of OGM-$q$} 
to \cblue{those of}
other choices of $\gamma=0,0.4,0.8$
\cblue{with $\beta=\betastar(\gamma)$ in~\eqref{eq:betaopt}.}
The \cblue{OGM-$q$} 
(\cblue{see} upper red curve in Fig.~\ref{fig:ogm,roots})
has 
\cblue{the largest value ($\rho^\star = \gamstar = 0.6$)
of $\rho(\T_{\lambda}(\Frac{1}{L},\beta,\gamma))$}
at both the smallest and the largest eigenvalues
\cblue{($\mu$ and $L$ respectively)},
unlike other choices of \gam (with $\betastar(\gamma)$)
where either $\rho(\T_{\mu}(\Frac{1}{L},\beta,\allowbreak\gamma))$
or $\rho(\T_L(\Frac{1}{L},\beta,\gamma))$
are the largest.
The other choices thus have a spectral radius
\cblue{$\rho(\T(\Frac{1}{L},\beta,\gamma))$}
larger than that of the \cblue{OGM-$q$}.

\cblue{The right plot of} Fig.~\ref{fig:ogm,roots}
illustrates \cblue{$\rho(\T_{\lambda}(\Frac{1}{L},\beta,\gamma))$ values}
for different choices of $\beta\cblue{=0,0.2,0.4,0.6}$
for given $\gamma = \gamstar$,
showing that suboptimal $\beta$ value
will slow down convergence,
\cblue{compared to the optimal $\beta^\star=0.4$}.
\cblue{AFM with
$(\alpha,\beta,\gamma) = (\Frac{1}{L},0,\gamma^\star)$
in Fig.~\ref{fig:ogm,roots}
is equivalent to 
AFM with $\paren{\frac{1}{L}(1 + \gamma^\star),0,0}$,
and this implies that AFM with $\beta=\gamma=0$ (\eg, GM)
may have} 
some modes for mid-valued $\lambda$ values
that will converge faster
than the accelerated methods,
\cblue{whereas} 
its overall convergence rate \cblue{(\ie, the spectral radius value)}
is worse.
Apparently no one method
can have superior convergence rates
for all modes. 

\cblue{Similarly,} although 
\cblue{OGM-$q$ has}
the smallest possible
spectral radius $\rho(\T(\cdot))$
\cblue{among known AFM},
the upper blue and red curves in
\cblue{the left plot of} Fig.~\ref{fig:ogm,roots},
\cblue{corresponding to FGM-$q$ and OGM-$q$ respectively,}
illustrate that 
\cblue{OGM-$q$}
will have modes
for large eigenvalues
that converge slower than with
\cblue{FGM-$q$}.
This behavior may be undesirable
when such modes \cblue{of large eigenvalues}
dominate the overall convergence behavior.

\cblue{The next section reveals}
that the convergence of the primary sequence $\{\y_k\}$ of 
\cblue{AFM with $\alpha=\Frac{1}{L}$}
is not governed by such modes \cblue{of large eigenvalues}
unlike \cblue{its} secondary sequence $\{\x_k\}$. 
\cblue{In addition,}
Fig.~\ref{fig:ogm,roots}
reveals change points across $\lambda$
meaning that there are different regimes;
\cblue{the next section} elaborates on this behavior,
building upon \cblue{the} dynamical system analysis of
\cblue{AFM with $\alpha=\Frac{1}{L}$ and $\gamma=0$ in}
\cite[Sec. 4]{odonoghue:15:arf}.

\subsection{Convergence Properties of \cblue{AFM}
When
$\alpha = \Frac{1}{L}$}
\label{sec:OGM,conv}

\cite[Sec. 4]{odonoghue:15:arf}
analyzed a constant-step \cblue{AFM with $\alpha=\Frac{1}{L}$ and $\gamma=0$} 
as a linear dynamical system
for minimizing a strongly convex quadratic function~\eqref{eq:quad},
and showed that there are three regimes of behavior for the system;
low momentum, optimal momentum, and high momentum regimes.
This section
similarly analyzes \cblue{AFM with $\alpha=\Frac{1}{L}$ and $\gamma\ge0$} 
to better understand its convergence behavior
when solving a strongly convex quadratic problem~\eqref{eq:quad},
complementing the previous
section's
spectral radius analysis of \cblue{AFM}.

We use the eigen-decomposition of $\Q = \V\Lam\V^\top$
with $\Lam := \diag{\lambda_i}$,
where the eigenvalues $\{\lambda_i\}$ are in an ascending order,
\ie, $\mu=\lambda_1\le\lambda_2\le\cdots\le\lambda_d=L$.
And for simplicity, we let $\vp = \Zero$
without loss of generality,
leading to $\x_* = \Zero$.
By defining
$\w_k := (w_{k,1},\cdots,w_{k,d})^\top = \V^\top \y_k \in \Reals^d$
and
$\vv_k := (v_{k,1},\cdots,v_{k,d})^\top = \V^\top \x_k \in \Reals^d$
as the mode coefficients
of the primary and secondary sequences respectively
and using~\eqref{eq:xii,update},
we have the following $d$ independently evolving
identical recurrence relations
for the evolution of $w_{\cdot,i}$ and $v_{\cdot,i}$ 
of the constant-step \cblue{AFM with $\alpha=\Frac{1}{L}$} respectively:
\begin{align}
w_{k+2,i} &= \paren{(1 + \beta)\paren{1 - \Frac{\lambda_i}{L}}
	- \gamma\Frac{\lambda_i}{L}}w_{k+1,i}
	- \beta\paren{1 - \Frac{\lambda_i}{L}}w_{k,i},
\label{eq:wrecur} \\
v_{k+2,i} &= \paren{(1 + \beta)\paren{1 - \Frac{\lambda_i}{L}}
	- \gamma\Frac{\lambda_i}{L}}v_{k+1,i}
	- \beta\paren{1 - \Frac{\lambda_i}{L}}v_{k,i}
\nonumber
,\end{align}
for $i=1,\ldots,d$,
although the initial conditions differ as follows:
\begin{align}
w_{1,i} = (1 - \Frac{\lambda_i}{L})w_{0,i},
\quad
v_{1,i} = ((1+\beta+\gamma)(1-\Frac{\lambda_i}{L}) - (\beta+\gamma))v_{0,i}
\label{eq:initial}
\end{align}
with $w_{0,i} = v_{0,i}$.
The convergence behavior of the $i$th mode of the dynamical system
of both $w_{\cdot,i}$ and $v_{\cdot,i}$
in \eqref{eq:wrecur}
is determined
by the characteristic polynomial~\eqref{eq:poly} 
with $\alpha = \Frac{1}{L}$
and $\lambda = \lambda_i$.
Unlike the previous sections
that studied only the worst-case convergence performance
using the largest absolute value of the roots
of the polynomial~\eqref{eq:poly},
we next discuss
the convergence behavior of \cblue{AFM}
more comprehensively
using~\eqref{eq:poly} with $\alpha=\Frac{1}{L}$ and $\lambda=\lambda_i$
for the two cases
1) $\lambda_i = L$ and 2) $\lambda_i < L$.

1) $\lambda_i = L$:
The characteristic polynomial~\eqref{eq:poly}
of the mode of $\lambda_i = L$
reduces to $r^2 + \gamma r=0$
with two roots $0$ and $-\gamma$
regardless of the choice of $\beta$.
Thus we have monotone
convergence
for this ($d$th) mode of the
dynamical system~\cite[Sec. 17.1]{chiang:84:fmo}:
\begin{align}
w_{k,d} = 0^k + c_d(-\gamma)^k,
\quad
v_{k,d} = 0^k + \hat{c}_d(-\gamma)^k
\label{eq:wkd}
,\end{align}
where $c_d$ and $\hat{c}_d$ are constants
depending on the initial conditions~\eqref{eq:initial}.
Substituting $w_{1,d} = 0$
and $v_{1,d} = -(\beta+\gamma)v_{0,d}$~\eqref{eq:initial}
into~\eqref{eq:wrecur}
yields 
\begin{align}
c_d = 0,
\quad
\hat{c}_d = v_{0,d}\paren{1 + \Frac{\beta}{\gamma}}
\label{eq:wkd,c}
,\end{align}
illustrating that the primary sequence $\{w_{k,d}\}$
reaches its optimum after one iteration,
whereas the secondary sequence $\{v_{k,d}\}$
has slow monotone convergence of the distance to the optimum,
while exhibiting undesirable oscillation
due to the term
$ (-\gamma)^k $,
corresponding
to overshooting over the optimum.

2) $\lambda_i < L$:
In \eqref{eq:ogmquadalp}
we found the optimal overall $\betastar$
for \cblue{AFM when $\alpha=\Frac{1}{L}$}.
One can alternatively explore
what the best value of $\beta$ would be
for any given mode of the system
for comparison.
The polynomial~\eqref{eq:poly} has repeated roots
for the following $\beta$,
corresponding to the smaller zero of
the discriminant
$\Delta(\Frac{1}{L},\beta,\gamma,\lambda_i)$
for given \gam and $\lambda_i$:
\begin{align}
\beta_i^\star(\gamma)
:= \Frac{\paren{1-\sqrt{(1+\gamma)\Frac{\lambda_i}{L}}}^2}{(1-\Frac{\lambda_i}{L})}
\label{eq:betaistar}
.\end{align}
This choice 
satisfies
$\betastar = \betastar(\gamstar)
= \beta_1^\star(\gamstar)$~\eqref{eq:ogmquadalp},
because $\lambda_1$ is the smallest eigenvalue.
Next we examine
the convergence behavior of \cblue{AFM with $\alpha=\Frac{1}{L}$ and $\gamma\ge0$}
in the following three regimes,
similar to \cblue{AFM with $\alpha=\Frac{1}{L}$ and $\gamma=0$}
in~\cite[Sec. 4.3]{odonoghue:15:arf}:\footnote{
For simplicity in the momentum analysis,
we considered values $\beta$
within $[0\;1]$,
containing the \cblue{standard} $\beta_k$ values 
in
\cblue{Tables~\ref{tab:alg1} and~\ref{tab:alg2}.}
\cblue{
This restriction  
excludes the effect of the $\beta$
that corresponds
to the larger zero of the discriminant
$\Delta(\Frac{1}{L},\beta,\gamma,\lambda_i)$ 
for given $\gamma$ and $\lambda_i$,
and that is larger than $1$.
Any $\beta$ greater than $1$ 
has $\rho(\T_{\lambda_i}(\Frac{1}{L},\allowbreak\beta,\gamma))$ values
(in~\eqref{eq:rho} with $\alpha=\Frac{1}{L}$)
that are
larger than those for $\beta\in[\beta_i^\star(\gamma)\;1]$
due to the quasi-convexity of 
$\rho(\T_{\lambda_i}(\Frac{1}{L},\allowbreak\beta,\gamma))$ over $\beta$.
}
}
\begin{itemize}[leftmargin=25pt]
\item $\beta<\beta_i^\star(\gamma)$: 
	low momentum, over-damped,
\item $\beta=\beta_i^\star(\gamma)$: 
	optimal momentum, critically damped,
\item $\beta>\beta_i^\star(\gamma)$: 
	high momentum, under-damped.
\end{itemize}

If $\beta \le \beta_i^\star(\gamma)$,
the polynomial~\eqref{eq:poly} has two real roots,
$r_{1,i}$ and $r_{2,i}$
where we omit $(\Frac{1}{L},\beta,\gamma,\lambda_i)$
in $r_{\cdot,i} = r_\cdot(\Frac{1}{L},\beta,\gamma,\lambda_i)$ for simplicity.
Then, the system evolves as~\cite[Sec. 17.1]{chiang:84:fmo}:
\begin{align}
w_{k,i} = c_{1,i}r_{1,i}^k + c_{2,i}r_{2,i}^k,
\quad
v_{k,i} = \hat{c}_{1,i}r_{1,i}^k + \hat{c}_{2,i}r_{2,i}^k
\label{eq:wki_over}
,\end{align}
where constants $c_{1,i}$, $c_{2,i}$, $\hat{c}_{1,i}$ and $\hat{c}_{2,i}$
depend on the initial conditions~\eqref{eq:initial}.
In particular, when $\beta = \beta_i^\star(\gamma)$ \eqref{eq:betaistar},
we have the repeated root:
\begin{align}
r_i^\star(\gamma) := 1 - \sqrt{(1+\gamma)\Frac{\lambda_i}{L}}
\label{eq:ris}
,\end{align}
corresponding to critical damping,
yielding the fastest monotone convergence
among \eqref{eq:wki_over} 
for any $\beta$ \st $\beta\le\beta_i^\star(\gamma)$.
This property is due to the quasi-convexity 
of $\rho(\T_{\lambda_i}(\Frac{1}{L},\allowbreak\beta,\gamma))$ over $\beta$.
If $\beta < \beta_i^\star(\gamma)$,
the system is over-damped, which corresponds to the low momentum regime,
where the system is dominated by the larger root
that is greater than $r_i^\star(\gamma)$ \cblue{\eref{eq:ris}},
and thus has slow monotone convergence.
However, depending on the initial conditions~\eqref{eq:initial},
the system may only be dominated by the smaller root,
as noticed for the case $\lambda_i=L$ in~\eqref{eq:wkd}
\cblue{and~\eqref{eq:wkd,c}}.
Also note that
the mode of $\lambda_i=L$ is always in the low momentum regime
regardless of the value of $\beta$.

If $\beta > \beta_i^\star(\gamma)$,
the system is under-damped,
which corresponds to the high momentum regime.
This means that the system evolves
as~\cite[Sec. 17.1]{chiang:84:fmo}:
\begin{align}
w_{k,i} = c_i\paren{\sqrt{\beta(1-\Frac{\lambda_i}{L})}}^k
		\cos(k\psi_i(\beta,\gamma)-\delta_i),
	\label{eq:wki_under} \\
v_{k,i} = \hat{c}_i\paren{\sqrt{\beta(1-\Frac{\lambda_i}{L})}}^k
	\cos(k\psi_i(\beta,\gamma)-\hat{\delta}_i),
\nonumber
\end{align}
where the
frequency of the oscillation
is given by
\begin{align}
\psi_i(\beta,\gamma)
	:= \cos^{-1}\paren{\Frac{\paren{(1+\beta)(1-\Frac{\lambda_i}{L}) - \gamma\Frac{\lambda_i}{L}}}
		{\paren{2\sqrt{\beta(1-\Frac{\lambda_i}{L})}}}}
\label{eq:psii}
,\end{align}
and $c_i$, $\delta_i$, $\hat{c}_i$ and $\hat{\delta}_i$
denote constants that depend on the initial conditions~\eqref{eq:initial};
in particular for $\beta\approx1$,
we have $\delta_i\approx0$ and $\hat{\delta}_i\approx0$
so we will ignore them.

Based on the above momentum analysis,
we categorize
the behavior of the $i$th mode
of \cblue{AFM} for each $\lambda_i$ \cblue{in Fig.~\ref{fig:ogm,roots}}.
Regimes with two curves and one curve 
\cblue{(over $\lambda$)}
in Fig.~\ref{fig:ogm,roots}
correspond to the low- and high-momentum regimes, respectively.
In particular, for $\beta = \betastar(\gamma)$
in \cblue{the left plot of} Fig.~\ref{fig:ogm,roots},
most $\lambda_i$ values 
\cblue{(satisfying $\beta>\beta_i^\star(\gamma)$)}
experience high momentum
(and the optimal momentum for 
$\lambda_i$ 
satisfying $\betastar(\gamma) = \beta_i^\star(\gamma)$,
\eg, $\lambda_i = \mu$),
whereas modes where $\lambda_i\approx L$
experience low momentum.
The fast convergence of
the primary sequence $\{w_{k,d}\}$
in \eref{eq:wkd} \cblue{and \eref{eq:wkd,c}}
generalizes to the case $\lambda_i \approx L$,
corresponding to the lower curves
in Fig.~\ref{fig:ogm,roots}.
In addition, for $\beta\cblue{=0,0.2}$ \cblue{that are}
smaller than $\betastar(\gamma)$
in \cblue{the right plot of} Fig.~\ref{fig:ogm,roots},
both $\lambda \approx \mu$ and $\lambda \approx L$
experience low momentum
so increasing $\beta$ improves the convergence rate.

Based on the quadratic analysis in this section,
we would like to use appropriately large $\beta$ and \gam
coefficients,
namely $(\betastar,\gamstar)$,
to have fast monotone convergence
(for the dominating modes).
However,
such values require knowing
the function parameter $q = \mu/L$
that is usually unavailable in practice.
Using OGM \cblue{(and OGM$'$) in Table~\ref{tab:alg1}}
without knowing $q$
will likely lead to oscillation
due to the high momentum (or under-damping)
for strongly convex functions.
The next section describes restarting schemes
inspired by~\cite{odonoghue:15:arf}
that we suggest to use with OGM
to avoid such oscillation
and thus heuristically accelerate
the rate of OGM
for a strongly convex quadratic function
and even
for a convex function that is locally strongly convex.

\mbox{}
\section{Restarting Schemes}
\label{sec:restart}

Restarting an algorithm 
(\ie, starting the algorithm again 
by using the current iterate as the new starting point)
after
a certain number of iterations
or when some restarting condition is satisfied
has been found useful,
\eg,
for the conjugate gradient method
\cite{powell:77:rpf,nocedal:06},
called ``fixed restart'' and ``adaptive restart'' respectively.
The fixed restart approach was also studied
for accelerated gradient schemes such as FGM
in~\cite{nesterov:13:gmf,nemirovski:94:emi}.
Recently
adaptive restart of FGM was shown
to provide dramatic practical acceleration
without requiring knowledge of function parameters%
~\cite{odonoghue:15:arf,giselsson:14:mar,su:16:ade}.
Building upon those ideas,
this section reviews and applies
restarting approaches for OGM.
A quadratic analysis in~\cite{odonoghue:15:arf}
justified using a restarting condition for FGM;
this section extends that analysis to OGM
by studying an observable quantity of oscillation
that serves as an indicator
for restarting the momentum of OGM.

\subsection{Fixed Restart}

Restarting an algorithm
every $k$ iterations
can yield a linear rate
for decreasing a function in \SL
\cite[Sec. 5.1]{nesterov:13:gmf}%
~\cite[Sec. 11.4]{nemirovski:94:emi}.
\cblue{Suppose one restarts OGM every $k$ (inner) iterations
by initializing the $(j+1)$th outer iteration using
$\x_{j+1,0} = \x_{j,k}$,
where
$\x_{j,i}$
denotes an iterate at the $j$th outer iteration
and $i$th inner iteration.}
Combining the OGM rate
\cblue{in Table~\ref{tab:alg1}}
and the strong convexity inequality~\eqref{eq:strcvx}
yields the following linear rate
\cblue{for each outer iteration of OGM with fixed restart}:
\begin{align}
f(\x_{j,k}) - f(\x_*)
 \le \frac{L ||\x_{j,0} - \x_*||^2}{k^2}
 \le \frac{2L}{\mu k^2} (f(\x_{j,0}) - f(\x_*))
\label{eq:fixed,cost}
.\end{align}
This \cblue{rate is faster} than the
\(
\Frac{4L}{\mu k^2}
\)
rate 
\cblue{of one outer iteration of FGM with fixed restart}
(using the FGM \cblue{rate in Table~\ref{tab:alg1}}).
\cblue{For a given $N=jk$ total number of steps,
a simple calculation
shows that
the} optimal restarting interval $k$
minimizing
\cblue{the rate $\paren{\Frac{2L}{(\mu k^2)}}^j$ after $N$ steps
(owing from~\eqref{eq:fixed,cost})
is}
$k_{\mathrm{fixed}}^{ } := e\sqrt{\Frac{2}{q}}$
\cblue{that does not depend on $N$,
where $e$ is Euler's number.}

There are two drawbacks of the fixed restart approach
\cite[Sec. 3.1]{odonoghue:15:arf}.
First,
computing the optimal interval
$k_{\mathrm{fixed}}^{ }$
requires knowledge of $q$
that is usually unavailable in practice.
Second, using a global parameter $q$
may be too conservative
when the iterates enter locally strongly convex region.
Therefore,
adaptive restarting~\cite{odonoghue:15:arf}
\cblue{is more} useful in practice,
which we review next and then apply to OGM.
The above two drawbacks
also apply to the \cblue{methods in Table~\ref{tab:coeff}}
that assume knowledge of the global parameter $q$.

\subsection{Adaptive Restart}

To circumvent the drawbacks of fixed restart,
\cite{odonoghue:15:arf}
proposes
the following two adaptive restart schemes
for FGM:
\begin{itemize}[leftmargin=25pt]
\item Function scheme for restarting (FR): restart whenever
\begin{align}
f(\y_{k+1}) > f(\y_k)
\label{eq:fr}
,\end{align}
\item Gradient scheme for restarting (GR): restart whenever
\begin{align}
\Inprod{-\nabla f(\x_k)}{\y_{k+1} - \y_k} < 0
\label{eq:gr}
.\end{align}
\end{itemize}
These schemes
heuristically improve convergence rates
of FGM 
and both performed similarly
well~\cite{odonoghue:15:arf,su:16:ade}.
Although
the function scheme
guarantees monotonic decreasing function values,
the gradient scheme has two advantages
over the function scheme~\cite{odonoghue:15:arf};
the gradient scheme involves only arithmetic operations
with already computed quantities,
and it is numerically more stable.

These two schemes encourage an algorithm to restart
whenever the iterates take a ``bad'' direction,
\ie,
when the function value increases
or the negative gradient and the momentum have an obtuse angle,
respectively.
However, a convergence proof that justifies
their empirical acceleration is yet unknown,
so~\cite{odonoghue:15:arf} 
analyzes such restarting schemes
for strongly convex quadratic functions.
An alternative scheme
in~\cite{su:16:ade} that restarts
whenever the magnitude of the momentum decreases, \ie,
$||\y_{k+1} - \y_k|| < ||\y_k - \y_{k-1}||$,
has a theoretical convergence analysis
for the function class \SL.
However,
empirically both the function and gradient schemes
performed better in~\cite{su:16:ade}.
Thus,
this paper
focuses on
adapting practical restart
schemes to OGM
and extending the
analysis in~\cite{odonoghue:15:arf} to OGM.
First we introduce a new additional adaptive scheme
designed specifically for 
\cblue{AFM with $\alpha=\Frac{1}{L}$ and $\gamma > 0$ (\eg, OGM).}

\subsection{Adaptive Decrease of \gam for 
\cblue{AFM with $\alpha=\Frac{1}{L}$ and $\gamma>0$}}
\label{sec:gr,gamma}

Sec.~\ref{sec:OGM,conv}
described
that the secondary sequence $\{\x_k\}$ of 
\cblue{AFM with $\alpha=\Frac{1}{L}$ and $\gamma>0$ (\eg, OGM)}
might experience overshoot and thus slow convergence,
unlike \cblue{its} primary sequence $\{\y_k\}$,
when the iterates enter 
a region where the mode of the largest eigenvalue dominates.
(Sec.~\ref{sec:result,case2} illustrates such an example.)
From~\eqref{eq:wkd},
the overshoot of $\x_k$
has magnitude proportional to $|\gamma|$,
yet a suitably large \gam,
such as \gamstar \eqref{eq:optogmalpbet},
is essential for overall acceleration.

To avoid (or reduce) such overshooting,
we suggest the following adaptive scheme:
\begin{itemize}[leftmargin=25pt]
\item Gradient scheme for decreasing \gam (GD\gam):
decrease \gam whenever
\begin{align}
\inprod{\nabla f(\x_k)}{\nabla f(\x_{k-1})} < 0
\label{eq:gr,gamma}
.\end{align}
\end{itemize}
Because
the primary sequence $\{\y_k\}$ of \cblue{AFM with $\alpha=\Frac{1}{L}$}
is unlikely to overshoot,
one could choose to simply
use the primary sequence $\{\y_k\}$
as algorithm output
instead of the secondary sequence $\{\x_k\}$.
However, if one needs
to use the secondary sequence of 
\cblue{AFM with $\alpha=\Frac{1}{L}$ and $\gamma>0$}
(\eg, Sec.~\ref{sec:nonstr}),
adaptive scheme
\eqref{eq:gr,gamma}
can help.

\subsection{Observable \cblue{AFM} Quantities \cblue{When $\alpha=\Frac{1}{L}$}}
\label{sec:OGM,obs}

This section revisits Sec.~\ref{sec:OGM,conv}
that suggested that
observing the evolution
of the mode coefficients $\{w_{k,i}\}$ and $\{v_{k,i}\}$
can help
identify the momentum regime.
However,
in practice
that evolution is unobservable
because the optimum $\x_*$ is unknown,
whereas
Sec.~\ref{sec:OGM,conv}
assumed $\x_* = \Zero$.
Instead we can observe the evolution of the function values,
which are related to the mode coefficients as follows:
\begin{align}
f(\y_k) = \frac{1}{2}\sum_{i=1}^d \lambda_i w_{k,i}^2,
\quad
f(\x_k) = \frac{1}{2}\sum_{i=1}^d \lambda_i v_{k,i}^2
\label{eq:cond,func}
,\end{align}
and also the inner products of the gradient and momentum, \ie,
\begin{align}
\inprod{-\nabla f(\x_k)}{\y_{k+1} - \y_k}
&= -\sum_{i=1}^d \lambda_i v_{k,i}(w_{k+1,i} - w_{k,i}), 
\label{eq:cond,grad1} \\
\inprod{\nabla f(\x_k)}{\nabla f(\x_{k-1})}
&= \sum_{i=1}^d \lambda_i^2 v_{k,i}v_{k-1,i}
\label{eq:cond,grad2}
.\end{align}
These quantities appear
in the conditions for the adaptive schemes
\eqref{eq:fr},~\eqref{eq:gr}, and~\eqref{eq:gr,gamma}.

One would like to increase
$\beta$ and \gam as much as possible for acceleration
up to \betastar and \gamstar~\eqref{eq:ogmquadalp}.
However, without knowing $q$ (and \betastar,\gamstar),
\cblue{using large $\beta$ and \gam} 
could end up placing
the majority of the modes
in the high momentum regime,
eventually leading to slow convergence
with oscillation as described in Sec.~\ref{sec:OGM,conv}.
To avoid such oscillation,
we
hope to detect it
using~\eqref{eq:cond,func} and~\eqref{eq:cond,grad1}
and
restart the algorithm.
We also hope to detect the overshoot~\eqref{eq:wkd} 
of the modes of the large eigenvalues (in the low momentum regime)
using~\eqref{eq:cond,grad2}
so that we can then decrease \gam and avoid such overshoot.

\cblue{The rest of this section}
focuses on
the case where
$\beta > \beta_1(\gamma)$
for given \gam,
when the most of the modes are in the high momentum regime.
Because the maximum
of $\rho(\T_\lambda(\Frac{1}{L},\beta,\gamma))$
occurs at the points $\lambda=\mu$ or $\lambda=L$,
we expect that%
~\eqref{eq:cond,func},~\eqref{eq:cond,grad1}, and~\eqref{eq:cond,grad2}
will be quickly dominated
by the mode of the smallest or the largest \cblue{eigen}values.
\cblue{Specifically,
plugging $w_{k,i}$ and $v_{k,i}$ 
in~\eqref{eq:wkd},~\eqref{eq:wkd,c} and~\eqref{eq:wki_under}
to~\eqref{eq:cond,func},~\eqref{eq:cond,grad1}, and~\eqref{eq:cond,grad2}
for only the (dominating) mode of the smallest and the largest eigenvalues
\cblue{($\lambda_1=\mu$ and $\lambda_d=L$ respectively)}
leads to
the following approximations:}
\begingroup
\allowdisplaybreaks
\begin{align}
f(\y_k) \approx&\, \frac{1}{2}\mu c_1^2 \, \beta^k \, (1-\Frac{\mu}{L})^k \, \cos^2(k\psi_1),
\label{eq:mode} \\
f(\x_k) \approx&\, \frac{1}{2}\mu\hat{c}_1^2 \, \beta^k \,
 (1-\Frac{\mu}{L})^k \, \cos^2(k\psi_1)
+ \frac{1}{2}L\hat{c}_d^2\gamma^{2k}
\nonumber \\
\inprod{-\nabla f(\x_k)}{\y_{k+1} - \y_k}
\approx&\, - \mu c_1\hat{c}_1 \, \beta^k \, (1-\Frac{\mu}{L})^k
	\cos(k\psi_1) \nonumber \\ 
	&\, \times \paren{ \sqrt{\beta(1-\Frac{\mu}{L})} \cos((k+1) \psi_1) - \cos(k\psi_1) },
\nonumber \\
\inprod{\nabla f(\x_k)}{\nabla f(\x_{k-1})}
	\approx&\, \mu^2 \hat{c}_1^2 \, \beta^{k-\frac{1}{2}} \, (1-\Frac{\mu}{L})^{k-\frac{1}{2}}
			\cos(k\psi_1) \, \cos((k-1)\psi_1) \nonumber \\
			&\,
			- L^2 \hat{c}_d^2 \, \gamma^{2k-1}
	\nonumber
,\end{align}
\endgroup
where
$\psi_1 = \psi_1(\beta,\gamma)$ \cblue{in~\eqref{eq:psii}}.
\cblue{Furthermore,}
it is likely that
these expressions
will be dominated
by the mode
of either the smallest or largest eigenvalues,
\cblue{so we}
next analyze each case separately.

\subsubsection{Case 1: the Mode of the Smallest Eigenvalue Dominates}
\label{sec:OGM,obs1}

When the mode of the smallest eigenvalue dominates,
we further approximate
\eqref{eq:mode} as
\begingroup
\allowdisplaybreaks
\begin{align}
&f(\y_k) \approx \frac{1}{2}\mu c_1^2 \, \beta^k \, (1-\Frac{\mu}{L})^k \, \cos^2(k\psi_1),
\quad	
f(\x_k) \approx \frac{1}{2}\mu \hat{c}_1^2 \, \beta^k \, (1-\Frac{\mu}{L})^k \, \cos^2(k\psi_1),
	\nonumber \\
&\inprod{-\nabla f(\x_k)}{\y_{k+1} - \y_k} \label{eq:mode,s} \\
	&\;\qquad\approx - \mu c_1\hat{c}_1 \, \beta^k \, (1-\Frac{\mu}{L})^k
		 \, \cos(k\psi_1) \, (\cos((k+1)\psi_1) - \cos(k\psi_1))
		\nonumber \\
	&\;\qquad= 2 \mu c_1\hat{c}_1 \, \beta^k \, (1-\Frac{\mu}{L})^k
		\, \cos(k\psi_1) \sin((k+\Frac{1}{2})\psi_1)
		\, \sin(\Frac{\psi_1}{2})
		\nonumber \\
	&\;\qquad\approx 2\mu c_1\hat{c}_1 \, \sin(\Frac{\psi_1}{2})
		\, \beta^k \, (1-\Frac{\mu}{L})^k \, \sin(2k\psi_1)
		\nonumber
,\end{align}
\endgroup
using simple trigonometric identities
and the approximations
$\sqrt{\beta(1-\Frac{\mu}{L})} \approx 1$
and
$\sin(k\psi_1) \approx \sin((k+\Frac{1}{2})\psi_1)$
\cblue{for small $\mu$ (leading to small $\psi_1$ in~\eqref{eq:psii})}.
The values~\eqref{eq:mode,s}
exhibit oscillations at a frequency proportional to
\(
\psi_1(\beta,\gamma)
\)
in
\eqref{eq:psii}.
This oscillation
can be
detected by the conditions~\eqref{eq:fr} and~\eqref{eq:gr}
and is useful
in detecting
the high momentum regime
where a restart
can help improve the convergence rate.

\subsubsection{Case 2: the Mode of the Largest Eigenvalue Dominates}
\label{sec:OGM,obs2}

Unlike the primary sequence $\{\y_k\}$ 
of \cblue{AFM with $\alpha=\Frac{1}{L}$ (\eg, OGM)},
convergence of
\cblue{its} secondary sequence $\{\x_k\}$
may be dominated by the mode of the largest eigenvalue
in~\eqref{eq:wkd} \cblue{and~\eqref{eq:wkd,c}}.
By further approximating~\eqref{eq:mode} for 
the case when the mode of the largest eigenvalue dominates,
the function value
\(
f(\x_k) \approx \frac{1}{2} L \hat{c}_d^2 \, \gamma^{2k}
\)
decreases slowly but monotonically,
whereas
$f(\y_k) \approx f(\x_*) = 0$
and
\(
\inprod{-\nabla f(\x_k)}{\y_{k+1} - \y_k}
	\approx 0 
.\)
Therefore,
neither restart condition~\eqref{eq:fr} or~\eqref{eq:gr}
can detect such non-oscillatory observable values,
even though the secondary mode $\{w_{k,d}\}$
of the largest eigenvalue is oscillating
(corresponding to overshooting over the optimum).
However,
the inner product of two sequential gradients:
\begin{align}
\inprod{\nabla f(\x_k)}{\nabla f(\x_{k-1})}
	\approx -L^2 \hat{c}_d^2 \, \gamma^{2k-1}
\end{align}
can detect the overshoot of the secondary sequence $\{\x_k\}$,
suggesting that the algorithm should adapt by decreasing
\gam
when condition~\eqref{eq:gr,gamma} holds.
Decreasing \gam too much
may slow down the overall convergence rate
when the mode of the smallest eigenvalue 
\cblue{is not negligible.}
Thus,
we use~\eqref{eq:gr,gamma}
only
when using
the secondary sequence \cblue{$\{\x_k\}$} 
as algorithm output
(\eg, Sec.~\ref{sec:nonstr}).

\section{Proposed Adaptive Schemes for OGM}
\label{sec:propOGM}

\subsection{Adaptive Scheme of OGM for Smooth and Strongly Convex Problems}
\label{sec:OGMrestart}

\algref{alg:OGMrestart}
illustrates
a new adaptive version of OGM\cblue{$'$}
\cblue{(rather than OGM)}\footnote{
\label{ft:ogmp}
OGM requires choosing the number of iterations $N$ in advance
for computing $\theta_N$ in Table~\ref{tab:alg1},
which seems incompatible with adaptive restarting schemes.
\cblue{In contrast,} the parameters $t_k$
in Table~\ref{tab:alg1} and \algref{alg:OGMrestart}
\cblue{are independent of $N$.}
The fact that $\theta_N$
is larger than $t_N$
at the last ($N$th) iteration
helps to dampen
(by reducing the values of $\beta$ and \gam)
the final update
to guarantee a faster \cblue{(optimal)}
worst-case rate \cblue{for the last secondary iterate $\x_N$}.
This property
was studied in~\cite{kim:17:otc}.
We could perform one last update using $\theta_N$
after a restart condition is satisfied,
but this
step appears unnecessary because
restarting already has the effect of
\cblue{dampening (reducing $\beta$ and $\gamma$)}.
Thus,
\cblue{\algref{alg:OGMrestart} uses OGM$'$ instead
that uses $t_k$
and that has a worst-case rate 
that is similar to that of OGM.}}
that is used in our numerical experiments in Sec.~\ref{sec:result}.
When a restart condition is satisfied in~\algref{alg:OGMrestart},
we reset $t_k = 1$
to discard
the previous momentum
that has a bad direction.
When the decreasing \gam condition
is satisfied in~\algref{alg:OGMrestart},
we decrease $\sigma$
to suppress undesirable overshoot
of the secondary sequence $\{\x_k\}$.
Although
the analysis in Sec.~\ref{sec:quadanal}
considered only strongly convex quadratic functions,
the numerical experiments in Sec.~\ref{sec:result}
illustrate that the adaptive scheme
is also useful
more generally
for smooth convex functions in \FL,
as described in~\cite[Sec. 4.6]{odonoghue:15:arf}.

\begin{algorithm}[htbp] 
\caption{OGM\cblue{$'$} with restarting momentum and decreasing \gam}
\label{alg:OGMrestart}
\begin{algorithmic}[1]
\State {\bf Input:}
	$f\in\SL$ or $\FL$,
	$\x_{-1} = \x_0 = \y_0\in\Reals^d$,
	$t_0=\sigma=1$, $\bsig \in [0,\;1]$.
\For{$k \ge 0$}
\State $\y_{k+1} = \x_k - \frac{1}{L}\nabla f(\x_k)$
\If{$f(\y_{k+1}) > f(\y_k)$
	(or $\Inprod{-\nabla f(\x_k)}{\y_{k+1} - \y_k} < 0$)}
	\Comment{Restart condition}
\State $t_k = 1$, $\sigma \leftarrow 1$
\ElsIf{$\Inprod{\nabla f(\x_k)}{\nabla f(\x_{k-1})} < 0$}
	\Comment{Decreasing \gam condition}
\State $\sigma \leftarrow \bsig \sigma$
\EndIf
\State $t_{k+1} = \frac{1}{2}\paren{1 + \sqrt{1 + 4t_k^2}}$
\State $\x_{k+1} = \y_{k+1}
	+ \frac{t_k-1}{t_{k+1}}(\y_{k+1} - \y_k)
	+ \sigma\frac{t_k}{t_{k+1}}(\y_{k+1} - \x_k)$
\EndFor
\end{algorithmic}
\end{algorithm}

\subsection{
Adaptive Scheme of a Proximal Version of OGM for Nonsmooth Composite Convex Problems
}
\label{sec:nonstr}

Modern applications often involve
nonsmooth composite convex problems:
\begin{align}
\argmin{\x} \; \{F(\x) := f(\x) + \phi(\x)\}
\label{eq:nonsmooth}
,\end{align}
where $f \in \FL$ is a smooth convex function
(typically not strongly convex)
and $\phi\in\Finf$ is a convex function
that is possibly nonsmooth
and ``proximal-friendly''~\cite{combettes:11:psm},
such as the $\ell_1$ regularizer
$\phi(\x) = ||\x||_1$.
Our numerical experiments in Sec.~\ref{sec:result}
show that
a new adaptive version
of a proximal variant of OGM
can be useful for solving such problems.

To solve~\eqref{eq:nonsmooth},
\cite{beck:09:afi} developed
a fast proximal gradient method,
popularized under the name
fast iterative shrinkage-thresh\allowbreak olding algorithm (FISTA).
\cblue{FISTA has the same rate as FGM in Table~\ref{tab:alg1}
for solving~\eqref{eq:nonsmooth},
by simply replacing the line 3 of \algref{alg:ogm}
with FGM coefficients
by $\y_{k+1} = \prox_{\alpha\phi}(\x_k - \alpha\nabla f(\x_k))$,
where the proximity operator is defined as
$
\prox_{h}(\z)
        := \argmin{\x}\{\frac{1}{2}||\z - \x||^2 + h(\x)\}
.$}
Variants of FISTA
with adaptive restart
were studied in~\cite[Sec. 5.2]{odonoghue:15:arf}.

Inspired by the fact that OGM
\cblue{has a worst-case rate} faster than FGM,
\cite{taylor:17:ewc}
studied a proximal variant\footnote{
Applying the proximity operator
to the primary sequence $\{\y_k\}$ of OGM,
similar to the extension of FGM to FISTA,
leads to
a poor worst-case rate~\cite{taylor:17:ewc}.
Therefore,~\cite{taylor:17:ewc}
applied the proximity operator
to the secondary sequence of OGM
and showed numerically
that this version has a worst-case rate
about twice faster than that of FISTA.
} of OGM (POGM).
It is natural
to pursue acceleration of POGM\footnote{
\cblue{Like OGM,} 
POGM
in~\cite[\cblue{Sec.~4.3}]{taylor:17:ewc}
requires choosing the number of iterations $N$ in advance
for computing $\theta_N$,
\cblue{and this is incompatible with adaptive restarting schemes.
Therefore, analogous to using OGM$'$ instead of OGM
for an adaptive scheme
in \algref{alg:OGMrestart} (see footnote~\ref{ft:ogmp}),
\algref{alg:POGMrestart}
uses a proximal version of OGM$'$
(rather than the POGM in~\cite{taylor:17:ewc})
with restart.
An extension of OGM$'$ (without restart)
to a proximal version
with a fast worst-case rate
is unknown yet}}
by using
variations of any (or all)
of the three adaptive schemes
\eqref{eq:fr},
\eqref{eq:gr},
\eqref{eq:gr,gamma},
as illustrated in \algref{alg:POGMrestart}.
\cblue{Regarding} a function restart condition for POGM,
we use
\begin{align}
F(\x_{k+1}) > F(\x_k)
\end{align}
instead of $F(\y_{k+1}) > F(\y_k)$,
because $F(\y_k)$ can be unbounded
(\eg, $\y_k$ can be unfeasible for constrained problems).
For gradient conditions of POGM,
we consider the composite gradient mapping
\cblue{$G(\x_k) \in \nabla f(\x_k) + \partial \phi(\x_{k+1})$} 
in \algref{alg:POGMrestart}
that differs from the standard composite gradient mapping
in~\cite{nesterov:13:gmf}.
We then use the gradient conditions
\begin{align}
\Inprod{-G(\x_k)}{\y_{k+1} - \y_k} < 0,
\quad
\Inprod{G(\x_k)}{G(\x_{k-1})} < 0
\label{eq:GR,gamma}
\end{align}
for restarting POGM or decreasing \gam of POGM
respectively.
Here POGM must output the secondary sequence $\{\x_k\}$
because the function value $F(\y_k)$ of the primary sequence
may be unbounded.
This situation
was the motivation
for~\eqref{eq:gr,gamma} 
\cblue{(and the second inequality of~\eqref{eq:GR,gamma})} 
and Sec.~\ref{sec:gr,gamma}.
When $\phi(\x) = 0$,
\algref{alg:POGMrestart} reduces to an algorithm that is similar to
\algref{alg:OGMrestart},
where only the location of the restart and decreasing $\gamma$ conditions differs.


\begin{algorithm}[htbp] 
\caption{\color{black} POGM\cblue{$'$} with restarting momentum and decreasing \gam}
\label{alg:POGMrestart}
\begin{algorithmic}[1]
\State {\bf Input:} $f\in\FL$, $\phi\in\Finf$, 
	$\x_{-1} = \x_0 = \y_0 = \u_0 = \z_0 \in\Reals^d$, \\
\qquad\quad $t_0=\zeta_0=\sigma=1$, $\bsig \in [0,\;1]$.
\For{$k \ge 0$}
\State $\u_{k+1} = \x_k - \frac{1}{L}\nabla f(\x_k)$
\State $t_{k+1} = \frac{1}{2}\paren{1 + \sqrt{1 + 4t_k^2}}$
\State $\z_{k+1} = \u_{k+1} + \frac{t_k-1}{t_{k+1}}(\u_{k+1} - \u_k)
	+ \sigma\frac{t_k}{t_{k+1}}(\u_{k+1} - \x_k)
	- \frac{t_k-1}{t_{k+1}}\frac{1}{L\zeta_k}(\x_k - \z_k)$
\State $\zeta_{k+1} = \frac{1}{L}
	\paren{1 + \frac{t_k-1}{t_{k+1}} + \sigma\frac{t_k}{t_{k+1}}}$
\State $\x_{k+1} = \prox_{\zeta_{k+1}\phi}(\z_{k+1})$
\State $G(\x_k) = \nabla f(\x_k) - \frac{1}{\zeta_{k+1}}(\x_{k+1} - \z_{k+1})$
\State $\y_{k+1} = \x_k - \frac{1}{L}G(\x_k)$
\If{$F(\x_{k+1}) > F(\x_k)$
        (or $\Inprod{-G(\x_k)}{\y_{k+1} - \y_k} < 0$)}
        \Comment{Restart condition}
\State $t_{k+1} = 1$, $\sigma \leftarrow 1$
\ElsIf{$\Inprod{G(\x_k)}{G(\x_{k-1})} < 0$}
        \Comment{Decreasing \gam condition}
\State $\sigma \leftarrow \bsig \sigma$
\EndIf
\EndFor
\end{algorithmic}
\end{algorithm}

\section{Numerical Results}
\label{sec:result}

This section shows the results of applying
OGM\cblue{$'$} and POGM\cblue{$'$}
with adaptive schemes
\cblue{in \algref{alg:OGMrestart} and \algref{alg:POGMrestart}}
to
various numerical examples including
both strongly convex quadratic problems
and non-strongly convex problems.\footnote{
Software for the algorithms and for producing the figures in Sec.~\ref{sec:result}
is available at 
\url{https://gitlab.eecs.umich.edu/michigan-fast-optimization/ogm-adaptive-restart}.
}
\cblue{(For simplicity, 
we omit the prime symbol of OGM$'$ and POGM$'$ with adaptive restart
hereafter.)}
The results illustrate that 
OGM (or POGM) with adaptive schemes
converges faster
than FGM (or FISTA) with adaptive restart.
The plots show the decrease of $F(\y_k)$ of the primary sequence
for FGM (FISTA) and OGM unless specified.
For POGM, we use the secondary sequence $\{\x_k\}$ as an output 
and plot $F(\x_k)$,
since $F(\y_k)$ can be unbounded.

\subsection{Strongly Convex Quadratic Examples}

This section considers two types
of strongly convex quadratic examples,
where the mode of either the smallest eigenvalue
or the largest eigenvalue dominates,
providing examples of the analysis
in Sec.~\ref{sec:OGM,obs1} and~\ref{sec:OGM,obs2}
respectively.

\subsubsection{Case 1: the Mode of the Smallest Eigenvalue Dominates}

\begin{figure}[htbp] 
\begin{center}
\includegraphics[clip,width=0.55\textwidth]{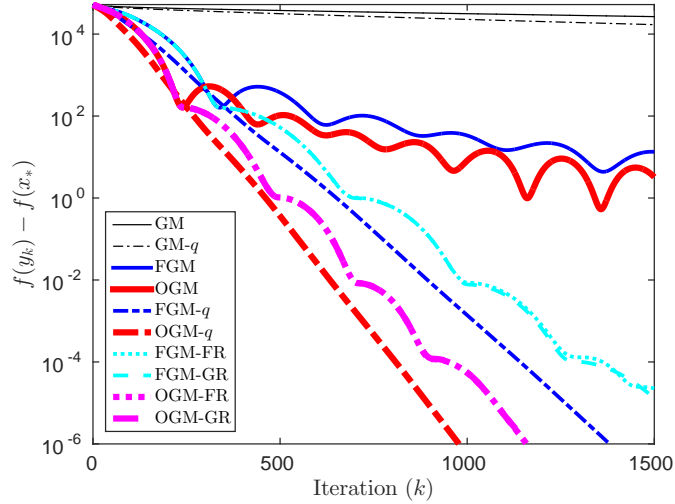}
\end{center}
\caption{Minimizing a strongly convex quadratic function -
Case 1: the mode of the smallest eigenvalue dominates.
(FGM-FR and FGM-GR are almost indistinguishable,
as are
OGM-FR and OGM-GR.)
}
\label{fig:toy}
\vspace{-3pt}
\end{figure}

Fig.~\ref{fig:toy} compares
GM, FGM and OGM,
with or without the knowledge of $q$,
for minimizing a strongly convex quadratic function~\eqref{eq:quad}
in $d=500$ dimensions with $q=10^{-4}$,
where we generated \A (for $\Q=\A\tr\A$) and \vp randomly.
As expected,
knowing $q$ accelerates convergence.

Fig.~\ref{fig:toy} also illustrates that
adaptive restart helps FGM and OGM
to nearly achieve
the fast linear converge rate
of their non-adaptive versions that know $q$.
As expected,
OGM \cblue{variants} converge faster than FGM \cblue{variants} for all cases.
In Fig.~\ref{fig:toy},
`FR' and `GR' stand for function restart~\eqref{eq:fr}
and gradient restart~\eqref{eq:gr},
respectively,
and both behave nearly the same.

\subsubsection{Case 2: the Mode of the Largest Eigenvalue Dominates}
\label{sec:result,case2}

Consider the strongly convex quadratic function
with
\(
\Q = 
\left[ \begin{array}{cc} q & 0 \\ 0 & 1 \end{array}\right]
,\)
$q = 0.01$,
$\vp = \Zero$
and $\x_* = \Zero$.
When starting the algorithm
from the initial point
\(
\x_0 = (0.2,\; 1)
,\)
the secondary sequence $\{\x_k\}$ of OGM-GR\footnote{
Fig.~\ref{fig:worst,toy}
only compares the results of the gradient restart (GR) scheme for simplicity,
where the function restart (FR) behaves similarly.
}
(or equivalently OGM-GR-GD\gam$(\bsig=1.0)$)
is dominated by the mode of largest eigenvalue in Fig.~\ref{fig:worst,toy},
illustrating
the analysis of Sec.~\ref{sec:OGM,obs2}.
Fig.~\ref{fig:worst,toy} illustrates that
the primary sequence of OGM-GR converges faster
than that of FGM-GR,
whereas the secondary sequence of OGM-GR
initially converges
even slower than GM.
To deal with such slow convergence
coming from the overshooting behavior
of the mode of the largest eigenvalue
of the secondary sequence of OGM,
we employ the decreasing \gam scheme in~\eqref{eq:gr,gamma}.
Fig.~\ref{fig:worst,toy} shows that using
$\bsig < 1$ in~\algref{alg:OGMrestart}
leads to overall faster convergence
of the secondary sequence $\{\x_k\}$
than the standard OGM-GR
where $\bsig = 1$.
We leave optimizing the choice of \bsig
or studying other strategies for decreasing \gam
as future work.

\begin{figure}[htbp] 
\begin{center}
\includegraphics[clip,width=0.55\textwidth]{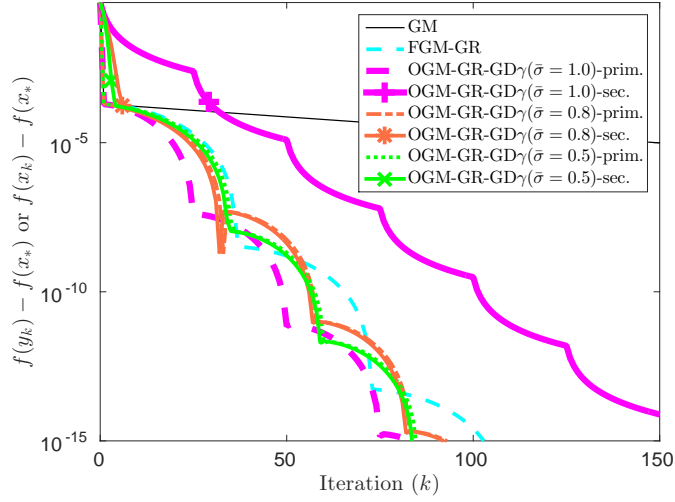}
\end{center}
\caption{Minimizing a strongly convex quadratic function -
Case 2: the mode of the largest eigenvalue dominates
for the secondary sequence $\{\x_k\}$ of OGM.
Using GD\gam~\eqref{eq:gr,gamma} with $\bsig < 1$ accelerates convergence
of the secondary sequence of OGM-GR,
where both the primary and secondary sequences
behave similarly
after first few iterations,
unlike $\bsig=1$.
}
\label{fig:worst,toy}
\vspace{-3pt}
\end{figure}

\subsection{Non-strongly Convex Examples}

This section applies adaptive OGM (or POGM)
to three non-strongly convex numerical examples
in~\cite{odonoghue:15:arf,su:16:ade}.
The numerical results show that
adaptive OGM (or POGM)
converges
faster than FGM (or FISTA) with adaptive restart.

\subsubsection{Log-Sum-Exp}

The following function
from~\cite{odonoghue:15:arf}
is smooth but non-strongly convex:
\[
f(\x) =
\eta\log\paren{\sum_{i=1}^m\exp\paren{\frac{1}{\eta}(\a_i^\top\x - b_i)}}
.\]
It approaches
$\max_{i = 1,\ldots,m} (\a_i^\top\x - b_i)$
as $\eta \to 0$.
Here,
$\eta$ controls the function smoothness
$L = \frac{1}{\eta}\lambda_{\max}(\A^\top\A)$
where $\A = [\a_1 \cdots \a_m]^\top \in \Reals^{m\times d}$.
The region around the optimum is approximately quadratic
since the function is smooth,
and thus the adaptive restart can be useful
without knowing the local condition number.

\begin{figure}[htbp]
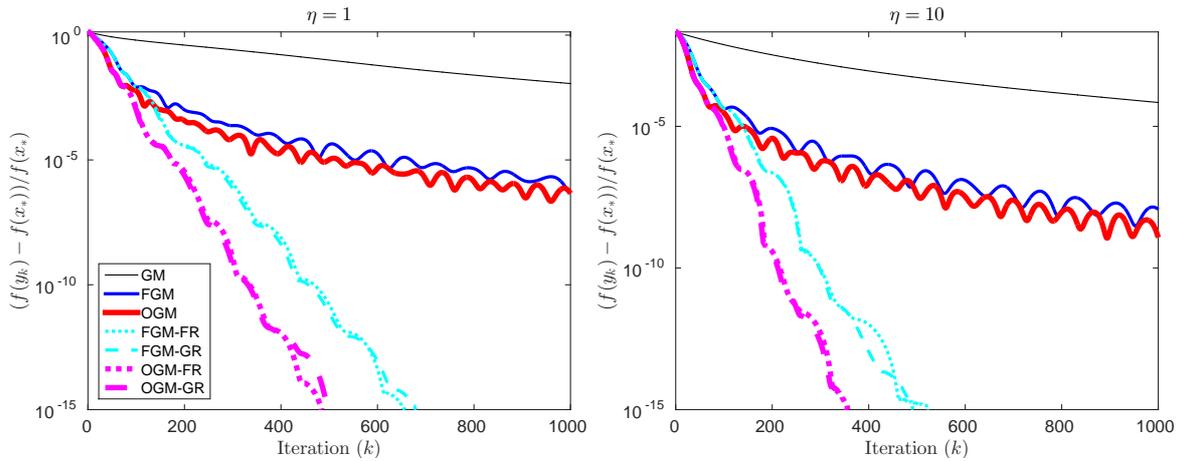
 
\begin{center}
\includegraphics[clip,width=0.48\textwidth]{fig,logsumexp1.eps}
\includegraphics[clip,width=0.48\textwidth]{fig,logsumexp2.eps}
\end{center}
\caption{Minimizing a smooth but non-strongly convex Log-Sum-Exp function.
}
\label{fig:logsumexp}
\vspace{-3pt}
\end{figure}

For $(m,d)=(100,20)$,
we randomly generated $\a_i\in\Reals^d$ and $b_i\in\Reals$
for $i=1,\ldots,m$,
and investigated $\eta=1,10$.
Fig.~\ref{fig:logsumexp}
shows that OGM with adaptive restart
converges faster
than FGM with the adaptive restart.
The benefit of adaptive restart is dramatic here;
apparently
FGM and OGM enter
a locally strongly convex region
after about $100-200$ iterations,
where adaptive restart
then provide a fast linear rate.

\subsubsection{Sparse Linear Regression}
\label{sec:sparselr}

Consider the following cost function
used for sparse linear regression:
\[
f(\x) = \frac{1}{2}||\A\x - \bb||_2^2
,\quad
\phi(\x) = \tau ||\x||_1
\label{eq:sparselr}
,\] 
for $\A\in\Reals^{m\times d}$,
where $L = \lambda_{\max}(\A^\top\A)$
and the parameter $\tau$ balances between
the measurement error and signal sparsity.
The proximity operator
becomes
a soft\hyp thresholding operator,
\eg,
$
\prox_{\zeta_{k+1}\phi}(\x)
	= \sgn(\x)\max\big\{|\x|-\zeta_{k+1}\tau,0\big\}
$.
The minimization seeks a sparse solution $\x_*$,
and often the cost function is strongly convex
with respect to the non-zero elements of $\x_*$.
Thus
we expect to benefit from adaptive restarting.

\begin{figure}[htbp]
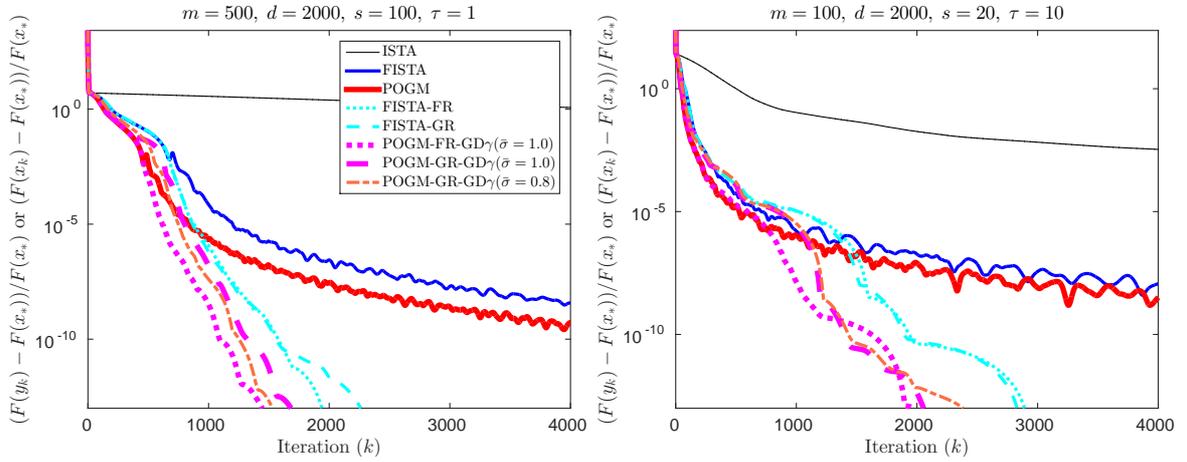
 
\begin{center}
\includegraphics[clip,width=0.48\textwidth]{fig,sparselr1.eps}
\includegraphics[clip,width=0.48\textwidth]{fig,sparselr2.eps}
\end{center}
\caption{Solving a sparse linear regression problem.
(ISTA is a proximal variant of GM.)}
\label{fig:sparselr}
\vspace{-3pt}
\end{figure}

For each choice of $(m, d, s, \tau)$ in Fig.~\ref{fig:sparselr},
we generated an $s$-sparse true vector $\x_{\mathrm{true}}$
by taking the $s$ largest entries of a randomly generated vector.
We then simulated
$\bb = \A\x_{\mathrm{true}} + \omegaa$,
where the entries of matrix \A and vector \omegaa
were sampled from a zero-mean normal distribution
with variances $1$ and $0.1$ respectively.
Fig.~\ref{fig:sparselr}
illustrates that POGM with adaptive schemes
provide acceleration over FISTA with adaptive restart.
While Sec.~\ref{sec:OGM,conv} discussed the undesirable overshooting behavior
that a secondary sequence of OGM (or POGM) may encounter,
these examples rarely encountered such behavior.
Therefore the choice of \bsig in the adaptive POGM
was not significant in this experiment,
unlike Sec.~\ref{sec:result,case2}.

\subsubsection{Constrained Quadratic Programming}

\begin{figure}[htbp]
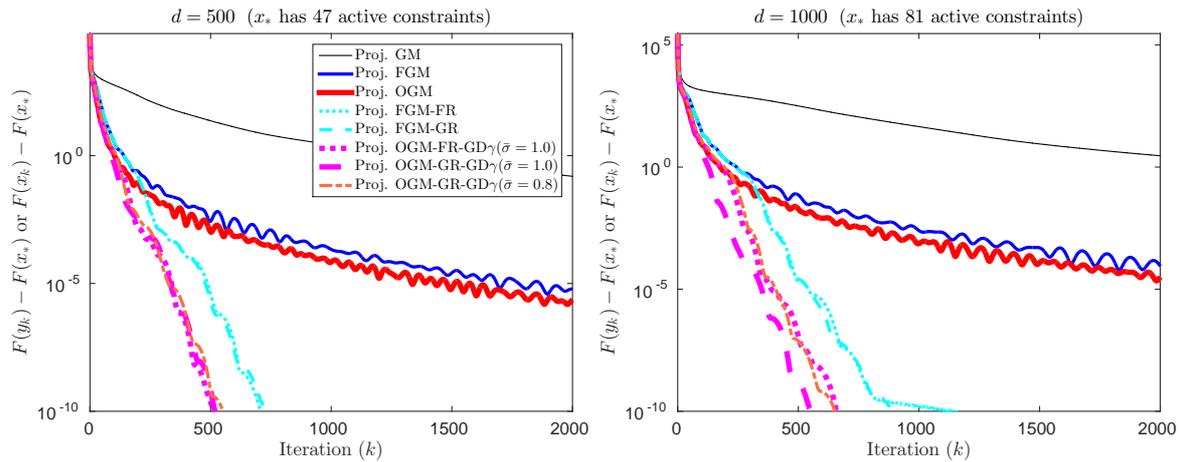
 
\begin{center}
\includegraphics[clip,width=0.48\textwidth]{fig,quad1.eps}
\includegraphics[clip,width=0.48\textwidth]{fig,quad2.eps}
\end{center}
\caption{Solving a box-constrained quadratic programming problem.}
\label{fig:quad}
\vspace{-15pt}
\end{figure}

Consider the following box-constrained quadratic program:
\[ 
f(\x) = \frac{1}{2}\x^\top\Q\x - \vp^\top \x,
\quad
\phi(\x) = \begin{cases}
                0, & \lll \preceq \x \preceq \u, \\
                \infty, & \text{otherwise},
        \end{cases}
,\] 
where $L = \lambda_{\max}(\Q)$.
The ISTA (a proximal variant of GM), 
FISTA and POGM use the projection operator:
\(
\prox_{\frac{1}{L}\phi}(\x)
	= \prox_{\zeta_{k+1}\phi}(\x)
	= \min\{\max\{\x, \lll\}, \u\}
.\)
Fig.~\ref{fig:quad} denotes
each algorithm by
a projected GM, a projected FGM, and a projected OGM respectively.
Similar to Sec.~\ref{sec:sparselr},
after the algorithm identifies the active constraints
the problem typically becomes a strongly convex quadratic problem
where we expect to benefit from adaptive restart.

Fig.~\ref{fig:quad} studies two examples
with problem dimensions $d = 500,1000$,
where we randomly generate a positive definite matrix \Q
having a condition number $10^7$ (\ie, $q = 10^{-7}$),
and a vector \vp.
Vectors \lll and \u
correspond to the interval constraints
$-1 \leq x_i \leq 1$ for $\x = \{x_i\}$.
The optimum $\x_*$ had $47$ and $81$ active constraints
out of $500$ and $1000$ respectively.
In Fig.~\ref{fig:quad},
the projected OGM with adaptive schemes
converged faster than
FGM with adaptive restart
and other non-adaptive algorithms.

\section{Conclusions}
\label{sec:conc}

We introduced adaptive restarting schemes
for the optimized gradient method (OGM)
to heuristically provide a fast linear convergence rate
when the function is strongly convex
or even when the function is not globally strongly convex.
The method resets
the momentum
when it makes a bad direction.
We provided a heuristic dynamical system analysis
to justify the practical acceleration of the adaptive scheme of OGM,
by extending the existing analysis of the fast gradient method (FGM).
On the way, we described 
\cblue{a new accelerated gradient method named OGM-$q$}
for strongly convex quadratic problems.
Numerical results illustrate that
the proposed adaptive approach
practically accelerates the convergence rate of OGM,
and in particular,
performs faster than FGM with adaptive restart.
An interesting open problem is to 
determine the worst-case rates
for OGM (and FGM) with adaptive restart.

\begin{acknowledgements}
This research was supported in part by NIH grant U01 EB018753.
\end{acknowledgements}

\bibliographystyle{spmpsci_unsrt} 
\bibliography{master}

\end{document}